\documentstyle{amsppt}
\pagewidth{5.59in} \pageheight{8.00in} \NoRunningHeads
\magnification = 1200

\topmatter
\title Eells-Sampson type theorems for subelliptic harmonic maps from sub-Riemannian manifolds*\endtitle
\author Yuxin Dong \endauthor
\thanks {*Supported by NSFC grant No. 11771087, and LMNS, Fudan.}
\endthanks
\abstract {In this paper, we consider critical maps of a
horizontal energy functional for maps from a sub-Riemannian
manifold to a Riemannian manifold. These critical maps are
referred to as subelliptic harmonic maps. In terms of the
subelliptic harmonic map heat flow, we investigate the existence
problem for subelliptic harmonic maps. Under the assumption that the
target Riemannian manifold has non-positive sectional curvature,
we prove some Eells-Sampson type existence results for this
flow when the source manifold is either a step-$2$ sub-Riemannian
manifold or a step-$r$ sub-Riemannian manifold whose
sub-Riemannian structure comes from a tense Riemannian foliation.
Finally, some Hartman type results are also established for the flow.}
\endabstract
\endsubjclass
\endkeywords
\endtopmatter
\document
\heading{\bf Introduction}
\endheading
\vskip 0.3 true cm

Sub-Riemannian geometry is a natural generalization of Riemannian
geometry, whose birth dates back to Carath\'{e}odory's 1909 seminal
paper on the foundations of Carnot thermodynamics. Geometric
analysis on sub-Riemannian manifolds has been received much
attention during the past decades (cf. [BBS1,2]). By a
sub-Riemannian manifold we mean a triple $(M,H,g_{H})$, where $M$
is a connected smooth manifold, $H$ is a subbundle of $TM$, and
$g_{H}$ is a smooth fiberwise metric on $H$. The subbundle $H$ is
usually assumed to have the bracket generating property for $TM$.
More precisely, one may introduce a generating order for the
sub-Riemannian manifold, that is, $M$ is called a step-$r$
sub-Riemannian manifold if sections of $H$ together with their Lie
brackets up to order $r$ spans $T_{x}M$ at each point $x$ (see \S
1 for the detailed definition). This is a remarkable property,
which makes both the geometry and analysis on sub-Riemannian
manifolds more interesting and rich.

The present paper is devoted to the study of a natural counterpart
of harmonic maps in the realm of sub-Riemannian geometry. A smooth
map $f:(M,H,g_{H})\rightarrow (N,h)$ from a sub-Riemannian
manifold with a smooth measure $d\mu$ to a Riemannian manifold is
called a subelliptic harmonic map if it is a critical map of the
following energy functional
$$
E_{H}(f)=\frac{1}{2}\int_{M}\mid df_{H}\mid ^{2}d\mu,\tag{0.1}
$$
where $df_{H}$ is the restriction of $df$ to $H$. To make the
above geometric variational problem manageable, we will restrict
our attention in this paper to a relative simple case that the
source sub-Riemannian manifold is endowed with a Riemannian
extension $g$ of $g_{H}$, and $d\mu=dv_{g}$ (the Riemannian volume
measure). We will find that the Euler-Lagrange-equations of the
functional (0.1) is a nonlinear subelliptic system of partial
differential equations (see \S 4 for its concrete expression)
$$
\tau _{H}(f)=0,  \tag{0.2}
$$
which justifies the terminology~for the critical map of $E_{H}$.
The principal part in (0.2) is actually the sub-Laplacian
$\bigtriangleup _{H}$, which is a hypoelliptic operator.

Recall that Jost-Xu [JX] first introduced subelliptic harmonic
maps associated with a H\"{o}rmander system of vector fields on a
domain of $R^{n} $ into Riemannian manifolds, and obtained an
existence and regularity theorem for these subelliptic maps under
Dirichlet condition and the same convexity condition of [HKW] on
the images. A related uniqueness result for subelliptic harmonic
maps in the sense of [JX] was given later by [Zh1]. As a global
formulation of Jost-Xu's subelliptic harmonic maps, E. Barletta at
al. introduced subelliptic harmonic maps from strictly
pseudoconvex CR manifolds into Riemannian manifolds, which were
referred to as pseudoharmonic maps in [BDU]; see also [DP] and
[Zh2] for some discussions on subelliptic harmonic maps from
almost contact Riemannian manifolds and sub-Riemannian manifolds
respectively. On the other hand, Wang [Wa] established some
regularity results for subelliptic harmonic maps from Carnot
groups, see also [HS], [ZF] for some regularity results of
subelliptic $p$-harmonic maps.

In the theory of harmonic maps, the Eells-Sampson theorem is a
fundamental theorem which has many essential applications in
Riemannian and K\"{a}hlerian geometry (cf.\ [JY], [Tol]). It
therefore seems natural and important to generalize this theorem
to the case of subelliptic harmonic maps from sub-Riemannian
manifolds. Note that step-$1$ sub-Riemannian manifolds are just
Riemannian manifolds. The simplest non-trivial sub-Riemannian
manifolds are step-$2$ sub-Riemannian manifolds, which includes
strictly pseudoconvex CR manifolds, contact metric manifolds,
quaternionic contact manifolds, or more general Heisenberg
manifolds, etc. (cf. [CC]). In [ChC], S. Chang and T. Chang gave
an Eells-Sampson type result for pseudoharmonic maps from compact
strictly pseudoconvex CR manifolds to compact Riemannian manifolds
with nonpositive curvature under an additional analytic condition
$[\bigtriangleup _{H},\xi ]=0$, where $\bigtriangleup _{H}$ and
$\xi $ are respectively the sub-Laplacian and Reeb vector field of
the source CR manifolds. Later, Y. Ren and G. Yang [RY] obtained a
general Eells-Sampson type result for pseduoharmonic maps without
Chang-Chang's condition. The main purpose in this paper is to
establish Eells-Sampson type theorems for subelliptic harmonic
maps from more general sub-Riemannian manifolds. Therefore we will
investigate the following subelliptic harmonic map heat flow
$$
\cases \frac{\partial f}{\partial t}=\tau_H(f)\\
 f|_{t=0}=\varphi\endcases\tag{0.3}
$$
for any given map $\varphi :(M,H,g_{H};g)\rightarrow (N,h)$. Our
main results include the short-time, long-time and homotopy
existence theorems for (0.3).

The paper is organized as follows. In \S 1 and \S 2, we collect
some basic notions and results about sub-Riemannian manifolds and
hypoelliptic PDEs from the literature. In \S 3, we first give the
structure equations of the generalized Bott connection
$\nabla^{\frak{B}}$ on a sub-Riemannian manifold $(M,H,g_{H};g)$;
and then introduce the second fundamental form of a map
$f:(M,H,g_{H};g)\rightarrow (N,h)$ with respect to the generalized
Bott connection on the source manifold and the Levi-Civita
connection on the target manifold. Using the moving frame method,
we are able to deduce some commutation relations for the
derivatives of the second fundamental form and thus some Bochner
type formulas for the map. In \S 4, we first give the
Euler-Lagrange-equations (0.2) in terms of the second fundamental
form of a map. Next, by means of the Nash embedding of the target
manifolds, we derive the explicit formulations for both (0.2) and
(0.3). \S 5 is devoted to existence problems. Using the heat
kernel associated with $\bigtriangleup _{H}-\partial _{t}$ and the
Duhamel's principle, we may establish a short time existence of
(0.2) for any initial map from a compact sub-Riemannian manifold
to a compact Riemannian manifold. When $N$ has nonpositive
curvature, we have the following long-time existence.

\proclaim{Theorem A} Let $(M,H,g_{H};g)$ be a compact
sub-Riemannian manifold and let $(N,h)$ be a compact Riemannian
manifold with nonpositive sectional curvature. Then for any smooth
map $\varphi :M\rightarrow N$, the subelliptic harmonic map heat
flow $(0.3)$ admits a global smooth solution $f:M\times \lbrack
0,\infty )\rightarrow N$.
\endproclaim

Under the nonpositive curvature condition on $N$, the above
theorem shows that the flow (0.3) does not blow up at any finite
time. Furthermore, in order to establish Eells-Sampson type
results for (0.3), one needs to have a uniform upper bound for the
energy density $e(f(\cdot ,t))$ of the solution $f(\cdot ,t)$ for
(0.3). We are able to give these uniform upper bounds in the
following two cases: the source manifolds are either step-$2$
sub-Riemannian manifolds or step-$r$ sub-Riemannian manifolds
whose sub-Riemannian structures come from some Riemannian
foliations. For both these cases, we have the Eells-Sampson type
results, which assert that there exists a sequence
$t_{i}\rightarrow \infty $ such that $f(x,t_{i}) \rightarrow
f_{\infty }(x)$ uniformly, as $t_{i}\rightarrow \infty $, to a $
C^{\infty }$ subelliptic harmonic map $f_{\infty }:M\rightarrow
N$. In \S 6, we establish Hartman type results for the subelliptic
harmonic map heat flow. Combining the Eells-Sampson and Hartman
type results, we have the following result for the first case.

\proclaim{Theorem B}Let $(M,H,g_{H};g)$ be a compact step-$2$
sub-Riemannian manifold and let $N$ be a compact Riemannian
manifold with non-positive sectional curvature. Then the
subelliptic harmonic map heat flow $(0.3)$ exists for all $t\in
\lbrack 0,\infty )$ and converges uniformly to a subelliptic
harmonic map $f_{\infty }$ as $t\rightarrow \infty $. In
particular, any map $\varphi \in C^{\infty }(M,N)$ is homotopic to
a $C^{\infty }$ subelliptic harmonic map.
\endproclaim

Riemannian foliations provide an important source of
sub-Riemannian manifolds. For a Riemannian foliation
$(M,g;\frak{F})$ with a bundle-like metric $g$, let
$H=(T\frak{F})^{\bot}$ (the horizontal subbundle of the foliation
$\frak{F}$ with respect to $g$) and $g_{H}$ be the restriction of
$g$ to $H$. Then we have a sub-Riemannian manifold $(M,H,g_{H};g)$
corresponding to $(M,g;\frak{F})$. The Riemannian foliation
$(M,g;\frak{F})$ will be said to be tense if the mean vector field
of $\frak{F}$ is parallel with respect to the Bott connection.
This is the second case in which we establish an Eells-Sampson
type result. Consequently we have

\proclaim{Theorem C}Let $(M,H,g_{H};g)$ be a compact
sub-Riemannian manifold corresponding to a tense Riemannian
foliation with the property that $H$ is bracket generating for
$TM$. Let $N$ be a compact Riemannian manifold with non-positive
sectional curvature. Then the subelliptic harmonic map heat flow
$(0.3)$ exists for all $t\in \lbrack 0,\infty )$ and converges
uniformly to a subelliptic harmonic map $f_{\infty }$ as
$t\rightarrow \infty $. In particular, any map $\varphi \in
C^{\infty }(M,N)$ is homotopic to a $C^{\infty }$ subelliptic
harmonic map.
\endproclaim

Hopefully these existence results will be useful for studying
either step-$2$ sub-Riemannian manifolds, such as contact and
quanternionic contact manifolds, or tense Riemannian foliations
with bracket generating horizontal subbundles. Besides their
possible geometric applications, we believe that it is reasonable
to investigate first the formulation for subellliptic harmonic
maps considered in this paper before studying more general
formulations, such as taking arbitrary smooth measures on the
source sub-Riemannian manifolds.

\heading{\bf 1. Sub-Riemannian geometry}
\endheading
\vskip 0.3 true cm

Let $M$ be a connected $(m+d)$-dimensional manifold of class
$C^{\infty }$ and let $H$ be a rank $m$ subbundle of the tangent
bundle $TM$. We say that $H$ satisfies the bracket generating
condition if vector fields which are sections of $H$ together with
all their brackets span $T_{x}M$ at each point $x$. More
precisely, for any $x\in M$ and any open neighborhood $U$ of $x$,
we let $\Gamma(U,H)$ denote the space of smooth sections of $H$ on
$U$, and define $\{\Gamma^{j}(U,H)\}_{j\geq 1}$ inductively by
$\Gamma^{j+1}(U,H)=\Gamma^{j}(U,H)+[\Gamma^{1}(U,H),\Gamma^{j}(U,H)]$
for each positive integer $j$, where
$\Gamma^{1}(U,H)=\Gamma(U,H)$. Here $[\cdot ,\cdot ]$ denotes the
Lie bracket of vector fields. By evaluating $\Gamma^{j}(U,H)$ at
$x$ , we have a subspace $H_{x}^{(j)}$ of the tangent space
$T_{x}M$, that is,
$$
H_{x}^{(j)}=\{X(x):X\in \Gamma^{j}(U,H)\}.  \tag{1.1}
$$
According to [St], [Mon], $H$ is said to be $r$-step bracket
generating for $TM$ if $H_{x}^{(r)}=T_{x}M$ for each $x\in M$.

A sub-Riemannian manifold is a triple $(M,H,g_{H})$, where $g_{H}$
is a fiberwise metric on the subbundle $H$. When $H$ is $1$-step
bracket generating, that is, $H=TM$, the sub-Riemannian manifold
is just a Riemannian manifold. Henceforth we will always assume
that $H$ satisfies the $r$-step bracket generating condition for
some $r\geq 2$. For a sub-Riemannian manifold, the subbundle $H$
is also referred to as a horizontal distribution. We say that a
Lipschitz curve $\gamma :[0,\delta ]\rightarrow M$ is horizontal
if $\gamma ^{\prime }(t)\in H_{\gamma (t)}$ a.e. in $[0,\delta ]$.
The sub-Riemannian metric $g_{H}$ induces a natural structure of
metric space, where the distance is the so-called
Carnot--Carath\'{e}odory distance
$$
\aligned d_{CC}(x_{0},x_{1})=\inf \{&\int_{0}^{\delta
}\sqrt{g_{H}(\gamma ^{\prime }(t),\gamma ^{\prime }(t))}dt\mid
\gamma :[0,\delta ]\rightarrow M\text{ is
a horizontal curve,} \\
&\gamma (0)=x_{0},\gamma (\delta )=x_{1}\}.
\endaligned\tag{1.2}
$$
By the theorem of Chow-Rashevsky ([Ch], [Ra]), there always exist
such curves joining $x_{0}$ and $x_{1}$, so the distance is finite
and continuous, and induces on $M$ the original topology. It turns
out that the distance $d_{CC}$ plays an essential role in
geometric analysis on sub-Riemannian manifolds. According to this
distance, we have a corresponding family of balls on $M$ given by
$$
B_{CC}(x,\delta )=\{y\in M\mid d_{CC}(x,y)<\delta \}.  \tag{1.3}
$$
These balls not only determine the metric topological properties
of $(M,d_{CC})$, but also reflect the non-isotropic feature of the
sub-Riemannian structure (cf. [NSW]).

One difficulty in sub-Riemannian geometry is the absence of a
canonical measure such as the Riemannian volume measure. Whenever
$M^{m+d}$ is endowed with a Riemannian metric $g$, we can compute
the volume of the $d_{CC}$-balls. One of the main results in [NSW]
is an estimate of the volume of these balls. To describe this
result, we choose a local frame field $\{e_{1},...,e_{m}\}$ of $H$
on a connected open subset $U\subset M$. Let
$$
\aligned
E^{(1)} &=\{e_{1},...,e_{m}\}, \\
E^{(2)} &=\{[e_{1},e_{2}],...,[e_{1},e_{m}],...,[e_{m-1},e_{m}]\},
\endaligned\tag{1.4}
$$
so that the components of $E^{(l)}$ are the commutators of length
$l$. Clearly $E^{(1)},...,E^{(l)}$ span $H^{(l)}$ at each point of
$U$ ($1\leq l\leq r$). Consequently, by the assumption for $H$, we
see that $E^{(1)},...,E^{(r)}$ span $TM$ at each point of $U$. Let
$Y_{1},...,Y_{q}$ be some enumeration of the components of
$E^{(1)},...,E^{(r)}$. A degree is assigned to each $Y_{i}$,
namely the corresponding length of the commutator. For each
$(m+d)$-tuple of integers $I=(i_{1},...,i_{m+d})$ with $1\leq
i_{j}\leq q$, following [NSW], one defines
$$
d(I)=\sum_{A=1}^{m+d}\deg (Y_{i_{A}})\qquad \text{and}\qquad
a_{I}(x)=|Y_{i_{1}}(x)\wedge \cdots \wedge Y_{i_{m+d}}(x)|_{g},\
x\in U. \tag{1.5}
$$
The Nagel-Stein-Wainger polynomial is defined by
$$
\Lambda (x,r)=\sum_{I}a_{I}(x)r^{d(I)},\quad r>0,  \tag{1.6}
$$
where the sum is over all $(m+d)$-tuples.

\proclaim{Theorem 1.1} (cf. [NSW]) Let $\{e_{i}\}_{i=1}^{m}$ be a
local frame field of $H$ on an open subset $U$ of $(M,g)$. Then,
for every open subset $V\ $of $U$ such that $\overline{V}\subset
U$ is compact, there exist constants $0<C,$ $R_{0}<1$, such
that for any $x\in V$, and $0<r\leq R_{0}$, one has%
$$
C\Lambda (x,r)\leq vol_{g}(B_{CC}(x,r))\leq C^{-1}\Lambda (x,r).
$$
\endproclaim

To describe the local growth order of $vol_{g}(B_{CC}(x,r))$, we
let
$$
Q(x)=\inf \{d(I)\mid a_{I}(x)\not=0\},\quad \text{ }Q=\sup
\{d(I)\mid |a_{I}(x)|\neq 0,\ x\in V\}\text{.}  \tag{1.7}
$$
According to [Ga], the numbers $Q(x)$ and $Q$ are respectively
called the pointwise homogeneous dimension of
$\{e_{i}\}_{i=1}^{m}$ at $x$ and the local homogeneous dimension
of $\{e_{i}\}_{i=1}^{m}$ on $U$. By the definitions of $Q(x)$ and
$Q$, one gets from (1.6) that
$$
t^{Q}\Lambda (x,r)\leq \Lambda (x,tr)\leq t^{Q(x)}\Lambda
(x,r),\quad 0<t\leq 1.  \tag{1.8}
$$

\proclaim{Corollary 1.2} (cf. also [Ga], [DGN]) For any $x\in V$,
$0<t\leq 1$, $0<r\leq R_{0}$, we have
$$
C_{1}t^{Q}\leq
\frac{vol_{g}(B_{CC}(x,tr))}{vol_{g}(B_{CC}(x,r))}\leq
C_{1}^{-1}t^{Q(x)}  \tag{1.9}
$$
where $C_{1}=C^{2}$. Besides, there exists a positive constant
$C_{2}$ such that
$$
vol_{g}(B_{CC}(x,r)\geq C_{2}r^{Q},\quad x\in V.  \tag{1.10}
$$
\endproclaim
\demo{Proof} Clearly (1.9) follows immediately from Theorem 1.1
and (1.8). Next Theorem 1.1 also yields
$$
vol_{g}(B_{CC}(x,r)\geq Cr^{Q}\sum_{I}a_{I}(x).
$$
Since $\sum_{I}a_{I}(x)>0$ on the compact set $\overline{V}$,
there exists a positive number $\widetilde{C}$ such that
$\sum_{I}a_{I}(x)\geq \widetilde{C} $ for any $x\in V$. Therefore
$vol_{g}(B_{CC}(x,r))\geq C_{2}r^{Q}$ with
$C_{2}=C\widetilde{C}$.\qed
\enddemo

For our purpose, we will consider compatible Riemannian metrics on
a sub-Riemannian manifold $(M,H,g_{H})$. A Riemannian metric $g$
on $M$ is called a Riemannian extension of $g_{H}$ if
$g|_{H}=g_{H}$. It is a known fact that such extensions always
exist (cf. [St]). Actually we may choose any Riemannian metric
$\widetilde{g}$ on $M$ and let $V$ be the orthogonal complement of
$H$ with respect to $\widetilde{g}$. Set $g_{V}=\widetilde{g}
|_{V}$. Then we have a Riemannian extension of $g_{H}$
$$
g=g_{H}+g_{V}  \tag{1.11}
$$
by requiring $g(u,v)=0$ for any $u\in H$ and $v\in V$. Clearly
such a Riemannian extension for $g_{H}$ is not unique. From now
on, we always fix a Riemannian extension $g$ on the sub-Riemannian
manifold $(M,H,g_{H})$, and consider the quadruple
$(M,H,g_{H};g)$. According to $g$, the tangent bundle $TM$ has the
following orthogonal decomposition:
$$
TM=H\oplus V.  \tag{1.12}
$$
The distribution $V$ will be referred to as the vertical
distribution or bundle on $(M,H,g_{H};g)$.

It would be convenient to introduce a suitable linear connection
compatible to the sub-Riemannian structure on $\left(
M,H,g_{H};g\right) $ in some sense. The generalized Bott
connection is one of such connections given by
$$
\nabla _X^{\frak{B}} Y=\cases &\pi _{H}(\nabla _X^R Y),\phantom{ab} X,Y\in \Gamma (H) \\
&\pi_H([X,Y]),\phantom{a} X\in \Gamma (V),Y\in \Gamma (H) \\
&\pi_V([X,Y]), \phantom{a} X\in \Gamma (H),Y\in \Gamma (V) \\
&\pi_V(\nabla_{X}^R Y),\phantom{a}  X,Y\in \Gamma (V)\endcases
\tag{1.13}
$$
where $\nabla^R$ denotes the Riemannian connection of $g$. Clearly
$\nabla^{\frak{B}}$ preserves the decomposition (1.12), and it
also satisfies
$$
\nabla _{X}^{\frak{B}}g_{H}=0\quad\text{ and }\quad\nabla
_{Y}^{\frak{B}}g_{V}=0 \tag{1.14}
$$
for any $X\in H$ and $Y\in V$. However, $\nabla^{\frak{B}}$ does
not preserve the Riemannian metric $g$ in general. The readers are
referred to [BF], [Ba2] for some discussions about this connection
on Riemannian foliations with totally geodesic leaves.

\example{Example 1.1} Let $G$ be a simply connected Lie group
whose Lie algebra $\frak{g}$ admits a direct sum decomposition of
vector spaces:
$$
\frak{g}=V_{1}\oplus V_{2}\oplus \cdots \oplus V_{r}\text{ \ \ \
}(r\geq 1)
$$
such that $[V_{1},V_{j}]=V_{j+1}$ \ for \ $1\leq j\leq r-1$ and $
[V_{1},V_{r}]=\{0\}$. Then $G$ is referred to as a Carnot group.
We may define a distribution $H$ on $G$ by
$H_{g}=dL_{g}(V_{1})\subset T_{g}G$, $\forall g\in G$. Let $g_{H}$
be a left-invariant metric on $H$. Clearly $(G,H,g_{H})$ is a
step-$r$ sub-Riemannian manifold. It is known that Carnot groups
play an important role in sub-Riemannian geometry and related
geometric analysis.
\endexample

\example{Example 1.2} Let $(M^{2n+1},\theta )$ be a (strict)
contact manifold, that is, $\theta $ is a global $1$-form
satisfying
$$
\theta \wedge (d\theta )^{n}\neq 0
$$
everywhere on $M$. Then the contact subbundle $H:=\ker \theta$ is
a $2$-step bracket generating subbundle of rank $2n$. The Reeb
vector field associated with $\theta $ is a unique vector field
$\xi $ on $M$ satisfying
$$
\theta (\xi )=1\text{ \ \ \ and \ \ \ }d\theta (\xi ,\cdot )=0.
$$
An almost complex structure $J$ in $H$ is said to be compatible
with $d\theta $ if
$$
d\theta (J\cdot ,J\cdot )=d\theta (\cdot ,\cdot )\text{ \ \ and \
\ }d\theta (J\cdot ,\cdot )>0.
$$
Then the contact subbundle $H$ and the Levi form $L_{\theta
}=d\theta (J\cdot ,\cdot )$ define a sub-Riemannian structure on
$M$.\ We extend $J$ to an endomorphism of $TM$ by setting $J\xi
=0$. The Webster metric defined by
$$
g_{\theta }=L_{\theta }+\theta \otimes \theta
$$
is a Riemannian extension of $L_{\theta }$. We call $(M,\theta
,\xi ,J,g_{\theta })$ a contact metric manifold. A contact metric
manifold $(M,\theta ,\xi ,J,g_{\theta })$ for which $J$ is
integrable is referred to as a strictly pseudoconvex CR manifold.
\endexample

\example{Example 1.3} (cf. [Biq1,2]) A quaternionic contact
manifold $M$ is a $(4n+3)$-dimensional manifold with a rank $4n$
distribution $H$ locally given as the kernel of 1-form $\eta
=(\eta ^{1},\eta ^{2},\eta ^{3})$ with values in $R^3$. In
addition, $H$ is equipped with a Riemannian metric $g_{H}$ and
three local almost complex structures $I_{i}$ ($i=1,2,3$)
satisfying the identities of the imaginary unit quaternions. These
structures also satisfy the following compatible conditions:
$g_{H}(I_{i}\cdot ,I_{i}\cdot )=g_{H}(\cdot ,\cdot )$ and $d\eta
_{i}=g(I_{i}\cdot ,\cdot )$. When the dimension of M is at least
eleven, Biquard [Biq1] also described the supplementary
distribution $V$ by the so-called Reeb vector fields $\{\xi
_{1},\xi _{2},\xi _{3}\}$. These Reeb vector fields are determined
by
$$
\eta _{s}(\xi _{k})=\delta _{sk},\quad (\xi _{s}\lrcorner d\eta
_{s})|_{H}=0,\quad (\xi _{s}\lrcorner d\eta _{k})|_{H}=-(\xi
_{k}\lrcorner d\eta _{s})|_{H}.
$$
Consequently $(H,g_{H})$ defines a $2$-step bracket generating
sub-Riemannian structure on $M$. Using the triple of Reeb vector
fields, we may extend $g_{H}$ to a Riemannian metric $g$ on $M$ by
requiring $span\{\xi _{1},\xi _{2},\xi _{3}\}=V\perp H$ and $g(\xi
_{s},\xi _{k})=\delta _{sk}$.
\endexample

\example{Example 1.4} (cf. [GW], [Mo]) A foliation on a manifold
is the collection of integral manifolds of an integrable
distribution on the manifold. Let $\frak{F}$ be a foliation on a
Riemannian manifold $(M,g)$. Set $V=T\frak{F}$, $H=V^{\bot }$
(w.r.t. $g$) and $g_{H}=g|_{H}$. Then $(H,g_{H})$ defines a
sub-Riemannian structure on $M$. The foliation is called a
Riemannian foliation if $\nabla _{\xi }^{\frak{B}}g_{H}=0$ for any
$\xi \in V$. In this case, following [Re], $g$ is referred to as a
bundle-like metric. Note that we are only interested in a
Riemannian foliation $\frak{F}$ whose horizontal distribution\ $H$
is bracket generating for $TM$ in this paper.
\endexample

For a sub-Riemannian manifold $(M,H,g_{H};g)$, we may define a
global vector field by
$$
\zeta =\pi _{H}(\sum_{\alpha }\nabla_{e_{\alpha }}^Re_{\alpha })
\tag{1.15}
$$
which will be called the mean curvature vector field of the
vertical distribution $V$. When $V$ is the tangent bundle of a
foliation $\frak{F}$ on $M$ as in Example 1.4, $\zeta $ is just
the usual mean curvature vector field along each leaf in
$(M,g,\frak{F})$. It is easy to verify by (1.13) that if
$\frak{F}$ is a Riemannian foliation with totally geodesic fibers,
then $\nabla ^{\frak{B}}$ is a metric connection for $g$ (cf.
[BF], [Ba2]).

\heading{\bf 2. Analysis for hypoelliptic operators}
\endheading
\vskip 0.3 true cm

In [H\"{o}], H\"{o}rmander considered the following type of
differential operator:
$$
\frak{D}=\sum_{i=1}^{m}X_{i}^{2}+Y  \tag{2.1}
$$
where $X_{1},...,X_{m},Y$ are smooth vector fields on a manifold
$\widetilde{M}$ with the property that their commutators up to
certain order span the tangent space at each point. He proved that
$\frak{D}$ is hypoelliptic in the sense that if $u$ is a
distribution defined on any open set $\Omega \subset
\widetilde{M}$, such that $\frak{D}u$ $\in C^{\infty }(\Omega )$,
then $u$ $\in C^{\infty }(\Omega )$. Due to this celebrated
result, hypoelliptic operators have since been the subject of
intense study (cf. [RS], [Br]). In following, we will discuss two
important hypoelliptic operators arising in sub-Riemannian
geometric analysis, namely, the sub-Laplacian and its heat
operator.

Let $(M^{m+d},H,g_{H};g)$ be a sub-Riemannian manifold with the
rank $m$ subbundle $H$ satisfying the $r$-step bracket generating
condition. A smooth vector field $X$ on $M$ is said to be
horizontal if $X_{p}\in H_{p}$ for each $p\in M$. For a smooth
function $u$, its horizontal gradient is the unique horizontal
vector field $\nabla ^{H}u$ satisfying $g_{H}(\nabla
^{H}u)_{q},X)=du(X)$ for any $X\in H_{q}$, $q\in M $. We choose a
local orthonormal frame field $\{e_{A}\}_{A=1}^{m+d}$ on an open
domain $\Omega $ of $(M,g)$ such that $span\{e_{i}\}_{i=1}^{m}=H$,
and thus $span\{e_{\alpha }\}_{\alpha =m+1}^{m+d}=V$. Such a frame
field is referred to as an adapted frame field for
$(M,H,g_{H};g)$. Consequently
$$
\nabla ^{H}u=\sum_{i=1}^{m}(e_{i}u)e_{i}.  \tag{2.2}
$$
Due to the H\"{o}rmander's condition, we see that $f$ is constant
if and only if $\nabla ^{H}u=0$.

By definition, in terms of the Riemannian connection $\nabla^R$,
the divergence of a vector field $X$ on $M$ is given by
$$
div_{g}X=\sum_{A=1}^{m+d}\{e_{A}\langle X,e_{A}\rangle -\langle
X,\nabla_{e_{A}}^Re_{A}\rangle \}. \tag{2.3}
$$
Then the sub-Laplacian of a function $u$ on $\left(
M,H,g_{H};g\right)$ is defined as
$$
\bigtriangleup_{H}u=div_{g}(\nabla ^{H}u).  \tag{2.4}
$$
Using the divergence theorem, we see that $\bigtriangleup _{H}$ is
a symmetric operator, that is,
$$
\int_{M}v(\bigtriangleup _{H}u)dv_{g}=\int_{M}u(\bigtriangleup
_{H}v)dv_{g}=-\int_{M}\mid \nabla ^{H}u\mid ^{2}dv_{g}  \tag{2.5}
$$
for any $u,v\in C_{0}^{\infty }(M)$. Using (2.2), (2.3), and
(1.15), we may rewrite (2.4) as
$$
\aligned \bigtriangleup_{H}u=&\sum_{i=1}^{m}\{e_{i}\langle \nabla
^{H}u,e_{i}\rangle -\langle \nabla ^{H}u,\nabla
_{e_{i}}^{\frak{B}}e_{i}\rangle \}-\langle
\nabla ^{H}u,\zeta \rangle \\
=&\sum_{i=1}^{m}e_{i}^{2}(u)-(\sum_{i=1}^{m}\nabla
_{e_{i}}^{\frak{B}} e_{i}+\zeta )u.\endaligned\tag{2.6}
$$
This shows that $\triangle _{H}$ is an operator of H\"{o}rmander
type, and thus it is hypoelliptic on $M$. Clearly the operator
$\bigtriangleup _{H}-\frac{\partial }{\partial t}$ is also an
operator given locally in the form of (2.1) with
$X_{1},...,X_{m},Y$ satisfying the H\"{o}rmander's condition on
$M\times R$. Therefore the heat operator corresponding to
$\bigtriangleup_{H}$ is hypoelliptic too.

In 1976, L. Rothschild and E.M. Stein [RS] established a more
precise regularity theory for hypoelliptic operators. Define
$$
S_{k}^{p}(\bigtriangleup _{H},\Omega)={\big \{}u\in L^{p}(\Omega
)\mid e_{i_{1}}\cdots e_{i_{s}}(u)\in L^{p}(\Omega ), \phantom{a}
1\leq i_{1},...,i_{s}\leq m,\phantom{a} 0\leq s\leq k {\big\}}
\tag{2.7}
$$
and
$$
\aligned S_{k}^{p}{\big(}\bigtriangleup _{H}-\frac{\partial
}{\partial t},\Omega \times (0,T){\big )}={\Big \{}&u\in
L^{p}{\big (}\Omega \times (0,T){\big )}\mid \partial
_{t}^{l}e_{i_{1}}\cdots e_{i_{s}}(u)\in L^{p}{\big (}\Omega \times
(0,T){\big )},\\ &1\leq i_{1},...,i_{s}\leq m,\ 2l+s\leq k {\Big
\}}\endaligned \tag{2.8}
$$
for any non-negative integer $k$. By the theory of Rothschild and
Stein, we have

\proclaim{Theorem 2.1} Let $\frak{D}=\bigtriangleup _{H}$ (resp.
$\bigtriangleup _{H}-\frac{
\partial }{\partial t}$) and $\widetilde{M}=\Omega $ (resp. $\Omega \times
(0,T)$). Suppose $f\in L_{loc}^{p}(\widetilde{M})$, and
$$
\frak{D}f=g\text{ \ \ \ on \ }\widetilde{M}.
$$
If $g\in S_{k}^{p}(\frak{D},\widetilde{M})$, then $\chi f\in
S_{k+2}^{p}(\frak{D},\widetilde{M})$ for any $\chi \in
C_{0}^{\infty }(\widetilde{M})$, $1<p<\infty $, $k=0,1,2,\cdots $.
In particular, the following inequality holds
$$
\parallel \chi f\parallel _{S_{k+2}^{p}(\frak{D},\widetilde{M})}\leq C_{\chi
}\left( \parallel g\parallel
_{S_{k}^{p}(\frak{D},\widetilde{M})}+\parallel f\parallel
_{L^{p}(\widetilde{M})}\right)
$$
where $C_{\chi }$ is a constant independent of $f$ and $g$.
\endproclaim

\remark{Remark 2.1} Let $L_{\alpha }^{p}(\widetilde{M})$,
$1<p<\infty $, be the classical Sobolev space. From [RS], we know
that $S_{k}^{p}(\bigtriangleup _{H}, \widetilde{M})\subset
L_{k/r}^{p}(\widetilde{M})$ for any $k\geq 0$, while $
S_{k}^{p}(\bigtriangleup _{H}-\frac{\partial }{\partial
t},\widetilde{M})\subset L_{k/r}^{p}(\widetilde{M})$ if $k$ is
even or a multiple of $r$. For any positive integer $l$, $\alpha
\in (0,1)$ and $1<p<\infty $, if $k$ is large enough, then $
S_{k}^{p}(\frak{D},\widetilde{M})\subset C^{l,\alpha
}(\widetilde{M})$ (the H\"{o}lder space) for
$\frak{D}=\bigtriangleup _{H}$ or $ \bigtriangleup
_{H}-\frac{\partial }{\partial t}$.
\endremark

Now we give some results about the heat kernel on compact
sub-Riemannian manifolds, which will be needed in \S 5. Let
$K(x,y,t)$ be the heat kernel for $\bigtriangleup _{H}$ on a
compact sub-Riemannian manifold $ (M,H,g_{H};g)$, that is,
$$
\cases &(\bigtriangleup _{H}-\frac{\partial }{\partial t})K(x,y,t)=0 \\
&\lim_{t\rightarrow 0}K(x,y,t)=\delta _{x}(y).
\endcases\tag{2.9}
$$
The readers may refer to [Ba1,3], [Bi] and [St] for the existence
of $K(x,y,t)$. We list some basic properties of $K(x,y,t)$ as
follows:
\newline\phantom{abcd}(1) $K(x,y,t)\in C^{\infty }(M\times M\times
R^{+})$;
\newline\phantom{abcd}(2) $K(x,y,t)=K(y,x,t)$ for $x,y\in M$ and
$t>0$;
\newline\phantom{abcd}(3) $K(x,y,t)>0$ for $x,y\in M$ and $t>0$;
\newline\phantom{abcd}(4) $\int_{M}K(x,y,t)dv_{g}(y)=1$ for any $x\in M$;
\newline\phantom{abcd}(5) $K(x,y,t+s)=\int_{M}K(x,z,t)K(y,z,s)dv_{g}(z)$ (semi-group
property).

The following result is a special case of a somewhat more general
theorem proved in [S\'{a}].

\proclaim{Theorem 2.2} (cf. [S\'{a}]) Let $K(x,y,t)$ be the heat
kernel of $\bigtriangleup _{H}$ on $(M,H,g_{H};g)$. Set
$w(x;\delta )=vol(B_{CC}(x;\delta ))$. Then
$$
|\nabla _{x}^{H}K(x,y,t)|\leq
A_{P}t^{-\frac{1}{2}}w(x;t^{1/2})^{-1}\left( 1+
\frac{d_{CC}(x,y)^{2}}{t}\right) ^{-P}  \tag{2.10}
$$
and
$$
K(x,y,t)\leq B_{P}w(x;t^{1/2})^{-1}\left(
1+\frac{d_{CC}(x,y)^{2}}{t}\right) ^{-P}  \tag{2.11}
$$
for $0<t<1$, all nonnegative integer $P$, and some positive
constants $A_{P}$ and $B_{P}$ depending on $P$, where $\nabla
_{x}^{H}$ denotes the horizontal gradient of $K$ with respect to
$x$.
\endproclaim

\proclaim{Lemma 2.3} For any $\beta \in (0,1/2)$, there exists a
$C_{\beta }>0$ such that
$$
\int_{0}^{t}\int_{M}|\nabla _{x}^{H}K(x,y,s)|dv_{g}(y)ds\leq
C_{\beta }t^{\beta }
$$
for $0<t<R_{0}$ for some positive constant $R_{0}$.
\endproclaim
\demo{Proof} Since $M$ is compact, there are two finite open
coverings $\{V_{a}\}$ and $ \{U_{a}\}$ ($a=1,...,l$) of $M$ such
that $\overline{V}_{a}\subset U_{a}$, $ \overline{V}_{a}$ is
compact, and Corollary 1.2 holds for each pair $ (V_{a},U_{a})$.
In particular, there exist positive constants $C_{a}$, $ D_{a} $
and $R_{a}$ such that for any $x\in V_{a}$, and $0<r<R_{a}$, one
has
$$
\frac{vol_{g}(B_{CC}(x,tr))}{vol_{g}(B_{CC}(x,r))}\geq
C_{a}t^{Q_{a}} \tag{2.12}
$$
for $0<t\leq 1$ and
$$
vol_{g}(B_{CC}(x,r))\geq D_{a}r^{Q_{a}}  \tag{2.13}
$$
where $Q_{a}$ is the local homogeneous dimension on $U_{a}$. For
any given $\beta \in (0,\frac{1}{2})$, we let $x\in V_{a}$ and
$\gamma _{a}=\frac{2\beta +Q_{a}-1}{2Q_{a}}$. \ Note that
$0<\gamma _{a}<\frac{1}{2}$, and thus $s^{(\frac{1}{2}-\gamma
_{a})}<1$ for any $0<s<1$. For any $0<s<R_{a}^{2}$, we obtain from
(2.12) that
$$
\aligned Vol_{g}(B_{CC}(x,s^{\frac{1}{2}}))
&=Vol_{g}(B_{CC}(x,s^{\frac{1}{2}-\gamma
_{a}}s^{\gamma _{a}})) \\
&\geq C_{a}s^{(\frac{1}{2}-\gamma
_{a})Q_a}vol_{g}(B_{CC}(x,s^{\gamma _{a}}))
\endaligned
$$
that is,
$$
\frac{vol_{g}(B_{CC}(x,s^{\gamma
_{a}}))}{Vol_{g}(B_{CC}(x,s^{\frac{1}{2}}))} \leq
C_{a}^{-1}s^{(\gamma _{a}-\frac{1}{2})Q_{a}}.  \tag{2.14}
$$
Taking a sufficiently large $P$ and using Theorem 2.2, (2.13) and
(2.14), we estimate the following integral for $0<t<R_{a}^{2}$:
$$
\aligned
&\int_{0}^{t}\int_{M}|\nabla _{x}^{H}K(x,y,s)|dv_{g}(y)ds \\
&\leq \int_{0}^{t}{\Big \{}\int_{B_{CC}(x,s^{\gamma
_{a}})}+\int_{M\,\diagdown \,B_{CC}(x,s^{\gamma _{a}})}{\Big
\}}|\nabla_{x}^{H}K(x,y,s)|dv_{g}(y)ds \\
&\leq A_{P}\int_{0}^{t}{\Big \{}\int_{B_{CC}(x,s^{\gamma
_{a}})}+\int_{M\,\smallsetminus \,B_{CC}(x,s^{\gamma _{a}})}{\Big
\}}\frac{(1+\frac{d_{CC}(x,y)^{2}}{s}
)^{-P}}{ s^{1/2}vol_{g}(B_{CC}(x,s^{1/2}))}dv_{g}(y)ds \\
&\leq A_{P}{\Big\{}\int_{0}^{t}\frac{vol_{g}(B_{CC}(x,s^{\gamma
_{a}}))}{s^{1/2}vol_{g}(B_{CC}(x,s^{1/2}))}ds+vol_{g}(M)\int_{0}^{t}\frac{
s^{P(1-2\gamma _{a})}}{s^{1/2}vol_{g}(B_{CC}(x,s^{1/2}))}ds{\Big\}} \\
&\leq A_{P}{\Big\{}\frac{1}{C_{a}}\int_{0}^{t}s^{(\gamma
_{a}-\frac{1}{2})Q_{a}-\frac{1}{2}}+\frac{vol_{g}(M)}{D_{a}}\int_{0}^{t}s^{P(1-2\gamma
_{a})-\frac{Q_{a}+1}{2}}ds{\Big\}} \\
&\leq \widetilde{C}_{a}t^{\beta }
\endaligned
$$
where $\widetilde{C}_{a}$ is a uniform positive constant. Set
$C_{\beta }=\max_{1\leq a\leq l}\{\widetilde{C}_{a}\}$ and
$R_{0}=\min_{1\leq a\leq l}\{R_{a}\}$. Then we complete the proof
of this lemma.\qed
\enddemo

In [Bo], Bony showed that the maximum principle holds for an
operator of H\"{o}rmander type. In following lemma, we provide
both a maximum principle (whose proof is routine), and a mean
value type inequality for subsolutions of the subelliptic heat
equation.

\proclaim{Lemma 2.4} Let $M$ be a compact sub-Riemannian manifold.
Suppose $\phi \ $is a subsolution of the subelliptic heat equation
satisfying
$$
\left( \bigtriangleup _{H}-\frac{\partial }{\partial t}\right)
\phi \geq 0
$$
on $M\times \lbrack 0,T)\ $with initial condition $\phi (x,0)=\phi
_{0}(x)$ for any $x\in M$. Then
$$
\sup_{M}\phi (x,t)\leq \sup_{M}\phi _{0}(x).
$$
Furthermore, if $\phi (x,t)$ is nonnegative, then there exist a
constant $B$ and an integer $Q$ such that
$$
\sup_{x\in M}\phi (x,t)\leq Bt^{-\frac{Q}{2}}\int_{M}\phi
_{0}(y)dv_{g}(y)
$$
for $0<t\leq \min \{R_{0}^{2},T\}$, where $R_{0}$ is as in Lemma
2.3.
\endproclaim
\demo{Proof} First we assume that $\phi $ is a subsolution of the
subelliptic heat equation. Set $c=\sup_{M}\phi _{0}(x)$. For any
fixed $\varepsilon >0$, one may introduce a function $\phi
_{\varepsilon }=\phi -\varepsilon (1+t)$. Clearly $\phi
_{\varepsilon }<c$ at $t=0$. We claim that $\phi _{\varepsilon
}<c$ for all $t>0$. In order to prove this, let us suppose the
result is false. This means that there exists $\varepsilon >0$
such that $\phi _{\varepsilon }\geq c$ somewhere in $M\times
\lbrack 0,T)$. Since $M$ is compact, there exists a point
$(x_{0},t_{0})\in M\times \lbrack 0,T)$ such that $\phi
_{\varepsilon }(x_{0},t_{0})=c$ and $\phi _{\varepsilon }(x,t)\leq
$ $c$ for all $x\in M$ and $t\in \lbrack 0,t_{0}]$. It follows
that $(\frac{\partial \phi _{\varepsilon }}{\partial
t})(x_{0},t_{0})\geq 0$ and $(\bigtriangleup _{H}\phi
_{\varepsilon })(x_{0},t_{0})\leq 0$, so that
$$
0>(\bigtriangleup _{H}\phi _{\varepsilon
})(x_{0},t_{0})-\varepsilon \geq (\frac{\partial \phi
_{\varepsilon }}{\partial t})(x_{0},t_{0})\geq 0
$$
which is a contradiction. Hence $\phi _{\varepsilon }<c$ on
$M\times \lbrack 0,T)$ for any $\varepsilon >0$. Since
$\varepsilon >0$ is arbitrary, we conclude that $\phi \leq c$ on
$M\times \lbrack 0,T)$. This proves the maximum principle.

Next we assume that $\phi $ is a nonnegative subsolution of the
subelliptic heat equation. Set
$$
\widetilde{\phi }(x,t)=\int_{M}K(x,y,t)\phi _{0}(y)dv_{g}(y).
\tag{2.15}
$$
Then $\widetilde{\phi }$ solves the subelliptic heat equation
$$
\left( \bigtriangleup _{H}-\frac{\partial }{\partial t}\right)
\widetilde{\phi }=0
$$
with initial data $\widetilde{\phi }(x,0)=\phi _{0}(x)$ for any
$x\in M$. By (2.15), we get
$$
\sup_{x\in M}\widetilde{\phi }(x,t)\leq \sup_{x,y\in
M}K(x,y,t)\int_{M}\phi _{0}(y)dv_{g}(y).  \tag{2.16}
$$
The semi-group property of $K(x,y,t)$ yields
$$
\aligned
K(x,y,t)&=\int_{M}K(x,z,\frac{t}{2})K(y,z,\frac{t}{2})dv_{g}(z) \\
&\leq \left( \int_{M}K^{2}(x,z,\frac{t}{2})dv_{g}(z)\right)
^{\frac{1}{2}
}\left( \int_{M}K^{2}(y,z,\frac{t}{2})dv_{g}(z)\right) ^{\frac{1}{2}} \\
&=K^{\frac{1}{2}}(x,x,t)K^{\frac{1}{2}}(y,y,t).
\endaligned\tag{2.17}
$$
According to Theorem 2.2, we have
$$
K(x,x,t)\leq \widetilde{B}\cdot vol_{g}{\big
(}B_{CC}(x,\sqrt{t}){\big )}^{-1} \tag{2.18}
$$
for some constant $\widetilde{B}$. Now we cover $M$ by two finite
open coverings $\{V_{a}\}_{a=1}^{l}$ and $\{U_{a}\}_{a=1}^{l}$ as
in the proof of Lemma 2.3. Let $Q_{a}$ be the local homogeneous
dimension on $U_{a}$. Set $Q=\max_{1\leq a\leq l}\{Q_{a}\}$. Then
we know from (2.13) that
$$
vol_{g}{\big (}B_{CC}(x,r){\big )}\geq Dr^{Q}  \tag{2.19}
$$
for $0<r\leq R_{0}=\min \{R_{a}\}$, where $D=\min_{1\leq a\leq
l}\{D_{a}\}$ and $Q=\max_{1\leq a\leq l}\{Q_{a}\}$. In terms of
(2.16), (2.17), (2.18) and (2.19), we conclude that
$$
\sup_{x\in M}\widetilde{\phi }(x,t)\leq
Bt^{-\frac{Q}{2}}\int_{M}\phi _{0}(y)dv_{g}(y)
$$
for $0<t\leq R_{0}^{2}$. Since $\phi $ is a subsolution, the
maximum principle implies that $\phi \leq \widetilde{\phi }$ for
$0<t<\min \{R_{0}^{2},T\}$. Hence we complete the proof of this
lemma.\qed
\enddemo

\heading{\bf 3. Second fundamental forms and their covariant
derivatives}
\endheading
\vskip 0.3 true cm

We will use the moving frame method to perform local computations
on maps from sub-Riemannian manifolds. For a sub-Riemannian
manifold $(M,H,g_{H};g)$, let us first give the structure
equations for the generalized Bott connection $\nabla ^{\frak{B}}$
defined by (1.13). Let $\{e_{A}\}_{A=1}^{m+p}$ be an adapted frame
field in $M$, and let $\{\omega ^{A}\}_{A=1}^{m+p}$ be its dual
frame field. From now on, we shall make use of the following
convention on the ranges of indices in $M$:
$$
\aligned
1\leq A, B, C&,...,\leq m+p;\quad 1\leq i, j, k,...,\leq m; \\
m&+1\leq\alpha, \beta, \gamma,...,\leq m+p,
\endaligned
$$
and we shall agree that repeated indices are summed over the
respective ranges. The connection $1$-forms $\{\omega_{A}^{B}\}$
of $\nabla ^{\frak{B}}$ with respect to $\{e_{A}\}_{A=1}^{m+d}$
are given by
$$
\nabla_{X}^{\frak{B}}e_{A}=\omega _{A}^{B}(X)e_{B} \tag{3.1}
$$
for any $X\in TM$. Since $\nabla^{\frak{B}}$ preserves the
decomposition (1.12), we have
$$
\nabla _{X}^{\frak{B}}e_{i}=\omega _{i}^{j}(X)e_{j},\quad \nabla
_{X}^{\frak{B}}e_{\alpha }=\omega _{\alpha }^{\beta }(X)e_{\beta }
\tag{3.2}
$$
and thus
$$
\omega _{i}^{\alpha }=0,\quad \omega _{\alpha }^{j}=0.  \tag{3.3}
$$
Let $T(\cdot,\cdot)$ and $R(\cdot,\cdot)$ be the torsion and
curvature of $\nabla^{\frak{B}}$ given respectively by
$$
\aligned
T(X,Y)&=\nabla _{X}^{\frak{B}}Y-\nabla _{Y}^{\frak{B}}X-[X,Y], \\
R(X,Y)Z&=\nabla _{X}^{\frak{B}}\nabla _{Y}^{\frak{B}}Z-\nabla
_{Y}^{\frak{B}}\nabla _{X}^{\frak{B}}Z-\nabla _{\lbrack
X,Y]}^{\frak{B}}Z
\endaligned\tag{3.4}
$$
where $X,Y,Z\in\Gamma(TM)$. Write
$$
\left. T(X,Y)=T^{A}(X,Y)e_{A},\quad R(X,Y)e_{A}=\Omega
_{A}^{B}(X,Y)e_{B}.\right.  \tag{3.5}
$$
Note that (3.2) implies
$$ \Omega _{i}^{\alpha }=\Omega _{\alpha
}^{j}=0.  \tag{3.6}
$$
As a linear connection, the structure equations of
$\nabla^{\frak{B}}$ are (cf. [KN])
$$
\aligned
d\omega ^{A}&=-\omega _{B}^{A}\wedge \omega ^{B}+T^{A}, \\
d\omega _{B}^{A}&=-\omega _{C}^{A}\wedge \omega _{B}^{C}+\Omega
_{B}^{A}.
\endaligned\tag{3.7}
$$

\proclaim{Lemma 3.1} For any $X,Y\in \Gamma(TM)$, we have
$$
T(X,Y)=-\pi _{V}([\pi _{H}(X),\pi _{H}(Y)])-\pi _{H}([\pi
_{V}(X),\pi _{V}(Y)]).
$$
\endproclaim
\demo{Proof} If $X,Y\in \Gamma (H)$, we verify by means of (1.13)
that
$$
\pi _{H}(T(X,Y))=\pi _{H}(\nabla _{X}^{R}Y-\nabla _{Y}^{
R}X-[X,Y])=0
$$
and
$$
\pi _{V}(T(X,Y))=-\pi _{V}([X,Y]).
$$
Similarly, if $X,Y\in \Gamma (V)$, then $\pi _{V}(T(X,Y))=0$ and
$\pi _{H}(T(X,Y))=-\pi _{H}\left( [X,Y]\right) $. Finally, if
$X\in \Gamma (V),Y\in \Gamma (H)$, then (1.13) implies directly
that $T(X,Y)=0$. Combining these cases, we prove this lemma.\qed
\enddemo

Using the dual frame field and Lemma 3.1, one may express the
torsion as
$$
\aligned T(\cdot,\cdot)=\frac{1}{2}( T_{ij}^{\alpha }&\omega
^{i}\wedge \omega ^{j}) \otimes e_{\alpha }+\frac{1}{2}( T_{\alpha
\beta }^{i}\omega ^{\alpha
}\wedge \omega ^{\beta }) \otimes e_{i} \\
&T_{ij}^{\alpha }=-T_{ji}^{\alpha },\qquad T_{\alpha \beta
}^{i}=-T_{\beta \alpha }^{i}.
\endaligned\tag{3.8}
$$
We also write
$$
\Omega _{B}^{A}=\frac{1}{2}R_{BCD}^{A}\omega ^{C}\wedge \omega
^{D},\quad R_{BCD}^{A}=-R_{BDC}^{A}.  \tag{3.9}
$$

Let $\left( N,h\right) $ be a Riemannian manifold and let
$\widetilde{\nabla }$ be its Riemannian connection. We choose an
orthonormal frame field $\{ \widetilde{e}_{I}\}_{I=1,...,n}$ in
$\left( N,h\right) $ and let $\{\widetilde{\omega }^{I}\}$ be its
dual frame field. The connection $1$-forms of $\widetilde{\nabla
}$ with respect to $\{\widetilde{e}_{I}\}_{I=1,...,n}$ are
$\{\widetilde{\omega }_{J}^{I}\}$. We will make use of the
following convention on the ranges of indices in $N$:
$$
I,J,K=1,...,n.
$$
The structure equations in $N$ are
$$
\aligned d\widetilde{\omega }^{I}&=-\widetilde{\omega
}_{K}^{I}\wedge \widetilde{
\omega }^{K} \\
d\widetilde{\omega }_{J}^{I}&=-\widetilde{\omega }_{K}^{I}\wedge
\widetilde{ \omega }_{J}^{K}+\widetilde{\Omega }_{J}^{I}
\endaligned\tag{3.10}
$$
where
$$
\widetilde{\Omega
}_{J}^{I}=\frac{1}{2}\widetilde{R}_{JKL}^{I}\widetilde{\omega
}^{K}\wedge \widetilde{\omega }^{L}.  \tag{3.11}
$$

For a smooth map $f:M\rightarrow N$, we have a connection $\nabla
^{\frak{B}}\otimes \widetilde{\nabla }^{f}$ in $T^{\ast }M\otimes
f^{-1}TN$, where $\widetilde{\nabla }^{f}$ denotes the pull-back
connection of $\widetilde{\nabla }$. Then the second fundamental
form with respect to the data $(\nabla
^{\frak{B}},\widetilde{\nabla }^{f})$ is defined by:
$$
\beta (f;\nabla ^{\frak{B}},\widetilde{\nabla
})(X,Y)=\widetilde{\nabla }_{Y}^{f}df(X)-df(\nabla
_{Y}^{\frak{B}}X). \tag{3.12}
$$
In terms of the frame fields in $M$ and $N$, the differential $df$
may be expressed as
$$
df=f_{A}^{I}\omega ^{A}\otimes \widetilde{e}_{I}.
$$
Consequently
$$
f^{\ast }\widetilde{\omega }^{I}=f_{A}^{I}\omega
^{A}=f_{i}^{I}\omega ^{i}+f_{\alpha }^{I}\omega ^{\alpha }.
\tag{3.13}
$$
By taking the exterior derivative of (3.13) and making use of the
structure equations in $M$ and $N$, we get
$$
Df_{A}^{I}\wedge \omega ^{A}+\frac{1}{2}f_{C}^{I}T_{AB}^{C}\omega
^{A}\wedge \omega ^{B}=0  \tag{3.14}
$$
where
$$
Df_{A}^{I}=df_{A}^{I}-f_{C}^{I}\omega
_{A}^{C}+f_{A}^{K}\widetilde{\omega }_{K}^{I}=f_{AB}^{I}\omega
^{B}.  \tag{3.15}
$$
Clearly the second fundamental form $\beta $ can be expressed as
$$
\beta =f_{AB}^{I}\omega ^{A}\otimes \omega ^{B}\otimes
\widetilde{e}_{I}. \tag{3.16}
$$
From (3.14), (3.15) and Lemma 3.1, it follows that
$$
\aligned
f_{ij}^{I}-f_{ji}^{I}&=f_{\alpha }^{I}T_{ij}^{\alpha } \\
f_{\alpha \beta }^{I}-f_{\beta \alpha }^{I}&=f_{k}^{I}T_{\alpha
\beta }^{k}
\\
f_{i\alpha }^{I}-f_{\alpha i}^{I}&=0.
\endaligned\tag{3.17}
$$
By taking the exterior derivative of (3.15), we deduce that
$$
\aligned Df_{AB}^{I}\wedge \omega ^{B}=&-f_{D}^{I}\Omega
_{A}^{D}+f_{A}^{K}\widetilde{\Omega }_{K}^{I}-f_{AD}^{I}T^{D} \\
=&-\frac{1}{2}f_{D}^{I}R_{ABC}^{D}\omega ^{B}\wedge \omega
^{C}+\frac{1}{2}
f_{A}^{K}\widetilde{R}_{KJL}^{I}f_{B}^{J}f_{C}^{L}\omega
^{B}\wedge \omega
^{C} \\
&-\frac{1}{2}f_{AD}^{I}T_{BC}^{D}\omega ^{B}\wedge \omega ^{C}
\endaligned\tag{3.18}
$$
where
$$
Df_{AB}^{I}=df_{AB}^{I}-f_{CB}^{I}\omega _{A}^{C}-f_{AC}^{I}\omega
_{B}^{C}+f_{AB}^{K}\widetilde{\omega }_{K}^{I}.  \tag{3.19}
$$
By putting
$$
Df_{AB}^{I}=f_{ABC}^{I}\omega ^{C},  \tag{3.20}
$$
we get from (3.18) the commutation relation
$$
f_{ABC}^{I}-f_{ACB}^{I}=f_{D}^{I}R_{ABC}^{D}+f_{AD}^{I}T_{BC}^{D}-f_{A}^{K}
\widetilde{R}_{KJL}^{I}f_{B}^{J}f_{C}^{L}.  \tag{3.21}
$$
For the map $f$, besides the differential $df$, one may also
introduce two partial differentials $df_{H}=df\mid_{H}\in \Gamma
(H^{\ast }\otimes f^{-1}TN)$ and $df_{V}=df\mid_{V}\in \Gamma
(V^{\ast }\otimes f^{-1}TN)$. By the definition of Hilbert-Schmidt
norm for a linear map, we have
$$
\mid df_{H}\mid ^{2}=\left( f_{i}^{I}\right) ^{2},\quad \mid
df_{V}\mid ^{2}=\left( f_{\alpha }^{I}\right) ^{2},\quad \mid
df\mid ^{2}=\left( f_{A}^{I}\right) ^{2}. \tag{3.22}
$$
Set
$$
e_{H}(f)=\frac{1}{2}\mid df_{H}\mid ^{2},\quad
e_{V}(f)=\frac{1}{2}\mid df_{V}\mid ^{2},\quad
e(f)=\frac{1}{2}\mid df\mid ^{2}.
$$
Now we want to derive the Bochner formulas of $\bigtriangleup
_{H}e_{H}(f)$, $\bigtriangleup _{H}e_{V}(f)$ and $\bigtriangleup
_{H}e(f)$. For a function $ u:M\rightarrow R$, one gets easily
from (2.6) and (3.12) that
$$
\aligned \bigtriangleup _{H}u&=\beta (u)(e_{k},e_{k})-\zeta (u)
\\& =u_{kk}-\zeta ^{k}u_{k}
\endaligned\tag{3.23}
$$
where $\zeta =\zeta ^{k}e_{k}$. Using (3.22), we compute
$$
e_{H}(f)_{k}=f_{i}^{I}f_{ik}^{I}  \tag{3.24}
$$
and
$$
{\big(}e_{H}(f){\big)}_{kk}=(f_{ik}^{I})^{2}+f_{i}^{I}f_{ikk}^{I}.
\tag{3.25}
$$
Consequently, in terms of (3.17) and (3.21), we derive that
$$
\aligned f_{ikk}^{I}&=[f_{ki}^{I}+f_{\alpha }^{I}T_{ik}^{\alpha
}]_{k}=f_{kik}^{I}+f_{\alpha k}^{I}T_{ik}^{\alpha }+f_{\alpha
}^{I}T_{ik,k}^{\alpha } \\
&=f_{kki}^{I}+f_{D}^{I}R_{kik}^{D}+f_{kD}^{I}T_{ik}^{D}-f_{k}^{K}\widetilde{R}
_{KJL}^{I}f_{i}^{J}f_{k}^{L}+f_{\alpha k}^{I}T_{ik}^{\alpha
}+f_{\alpha
}^{I}T_{ik,k}^{\alpha } \\
&=[f_{kk}^{I}-\zeta ^{k}f_{k}^{I}]_{i}+(\zeta
^{k}f_{k}^{I})_{i}+f_{D}^{I}R_{kik}^{D}+f_{kD}^{I}T_{ik}^{D}-f_{k}^{K}
\widetilde{R}_{KJL}^{I}f_{i}^{J}f_{k}^{L}+f_{\alpha
k}^{I}T_{ik}^{\alpha
}+f_{\alpha }^{I}T_{ik,k}^{\alpha } \\
&=\tau _{H,i}^{I}+\zeta _{,i}^{k}f_{k}^{I}+\zeta
^{k}f_{ki}^{I}+f_{D}^{I}R_{kik}^{D}+f_{kD}^{I}T_{ik}^{D}-f_{k}^{K}\widetilde{
R}_{KJL}^{I}f_{i}^{J}f_{k}^{L}+f_{\alpha k}^{I}T_{ik}^{\alpha
}+f_{\alpha }^{I}T_{ik,k}^{\alpha }
\endaligned\tag{3.26}
$$
where $\tau _{H}^{I}=f_{kk}^{I}-\zeta ^{k}f_{k}^{I}$ (see
Proposition 4.1 below for its geometric meaning). Then it follows
from (3.23), (3.25), (3.26) and (3.17) that
$$
\aligned
\bigtriangleup _{H}e_{H}(f)=&[e_{H}(f)]_{kk}-\zeta ^{k}[e_{H}(f)]_{k} \\
=&(f_{ik}^{I})^{2}+f_{i}^{I}f_{ikk}^{I}-\zeta ^{k}f_{i}^{I}f_{ik}^{I} \\
=&(f_{ik}^{I})^{2}+f_{i}^{I}\tau _{H,i}^{I}+f_{i}^{I}\zeta
_{,i}^{k}f_{k}^{I}+\zeta ^{k}f_{i}^{I}(f_{ki}^{I}-f_{ik}^{I}) \\
&+f_{i}^{I}f_{D}^{I}R_{kik}^{D}+f_{i}^{I}f_{kD}^{I}T_{ik}^{D}-f_{i}^{I}f_{k}^{K}%
\widetilde{R}_{KJL}^{I}f_{i}^{J}f_{k}^{L}+f_{i}^{I}f_{\alpha
k}^{I}T_{ik}^{\alpha }+f_{i}^{I}f_{\alpha }^{I}T_{ik,k}^{\alpha } \\
=&(f_{ik}^{I})^{2}+f_{i}^{I}\tau _{H,i}^{I}+f_{i}^{I}\zeta
_{,i}^{k}f_{k}^{I}+\zeta ^{k}f_{i}^{I}f_{\alpha
}^{I}T_{ki}^{\alpha
}+f_{i}^{I}f_{j}^{I}R_{kik}^{j} \\
&+f_{i}^{I}f_{k\alpha }^{I}T_{ik}^{\alpha
}-f_{i}^{I}f_{k}^{K}\widetilde{R}
_{KJL}^{I}f_{i}^{J}f_{k}^{L}+f_{i}^{I}f_{\alpha
k}^{I}T_{ik}^{\alpha
}+f_{i}^{I}f_{\alpha }^{I}T_{ik,k}^{\alpha } \\
=&(f_{ik}^{I})^{2}+f_{i}^{I}\tau _{H,i}^{I}+f_{i}^{I}\zeta
_{,i}^{k}f_{k}^{I}+\zeta ^{k}f_{i}^{I}f_{\alpha
}^{I}T_{ki}^{\alpha
}+f_{i}^{I}f_{j}^{I}R_{kik}^{j} \\
&+2f_{i}^{I}f_{\alpha k}^{I}T_{ik}^{\alpha
}-f_{i}^{I}f_{k}^{K}\widetilde{R}
_{KJL}^{I}f_{i}^{J}f_{k}^{L}+f_{i}^{I}f_{\alpha
}^{I}T_{ik,k}^{\alpha }.
\endaligned \tag{3.27}
$$
Similarly, using (3.17) and (3.21), we have
$$
\aligned {\big(} e_{V}(f){\big)}_{kk}=&[f_{\alpha }^{I}f_{\alpha
k}^{I}]_{k}=(f_{\alpha
k}^{I})^{2}+f_{\alpha }^{I}f_{\alpha kk}^{I} \\
=&(f_{\alpha k}^{I})^{2}+f_{\alpha }^{I}(f_{kk\alpha
}^{I}+f_{D}^{I}R_{k\alpha k}^{D}+f_{kD}^{I}T_{\alpha
k}^{D}-f_{k}^{K} \widetilde{R}_{KJL}^{I}f_{\alpha }^{J}f_{k}^{L}) \\
=&(f_{\alpha k}^{I})^{2}+f_{\alpha }^{I}\tau _{H,\alpha
}^{I}++f_{\alpha }^{I}\zeta _{,\alpha }^{k}f_{k}^{I}+f_{\alpha
}^{I}\zeta ^{k}f_{k\alpha
}^{I}+f_{\alpha }^{I}f_{D}^{I}R_{k\alpha k}^{D} \\
&+f_{\alpha }^{I}f_{kD}^{I}T_{\alpha k}^{D}-f_{\alpha
}^{I}f_{k}^{K}\widetilde{R}_{KJL}^{I}f_{\alpha }^{J}f_{k}^{L}.
\endaligned\tag{3.28}
$$
It follows that
$$
\aligned \bigtriangleup
_{H}[e_{V}(f)]&={\big(}e_{V}(f){\big)}_{kk}-\zeta ^{k}f_{\alpha
}^{I}f_{\alpha k}^{I} \\
&=(f_{\alpha k}^{I})^{2}+f_{\alpha }^{I}\tau _{H,\alpha
}^{I}+f_{\alpha }^{I}\zeta _{,\alpha }^{k}f_{k}^{I}+f_{\alpha
}^{I}f_{j}^{I}R_{k\alpha k}^{j}-f_{\alpha
}^{I}f_{k}^{K}\widetilde{R}_{KJL}^{I}f_{\alpha }^{J}f_{k}^{L}.
\endaligned\tag{3.29}
$$
From (3.27), (3.29), we conclude that
$$
\aligned \bigtriangleup _{H}e(f)=&(f_{ik}^{I})^{2}+(f_{\alpha
k}^{I})^{2}+f_{i}^{I}\tau _{H,i}^{I}+f_{\alpha }^{I}\tau
_{H,\alpha }^{I}+f_{i}^{I}\zeta _{,i}^{k}f_{k}^{I}+\zeta
^{k}f_{i}^{I}f_{\alpha
}^{I}T_{ki}^{\alpha } \\
&+f_{i}^{I}f_{j}^{I}R_{kik}^{j}+2f_{i}^{I}f_{\alpha
k}^{I}T_{ik}^{\alpha
}-f_{i}^{I}f_{k}^{K}\widetilde{R}_{KJL}^{I}f_{i}^{J}f_{k}^{L}+f_{i}^{I}f_{
\alpha }^{I}T_{ik,k}^{\alpha } \\
&+f_{\alpha }^{I}\zeta _{,\alpha }^{k}f_{k}^{I}+f_{\alpha
}^{I}f_{j}^{I}R_{k\alpha k}^{j}-f_{\alpha
}^{I}f_{k}^{K}\widetilde{R} _{KJL}^{I}f_{\alpha }^{J}f_{k}^{L}
\endaligned \tag{3.30}
$$
\proclaim{Lemma 3.2} Let $(M,H,g_{H};g)$ be a compact
sub-Riemannian manifold and let $(N,h)$ be a Riemannian manifold
with non-positive sectional curvature. Let $f:M\rightarrow N$ be a
smooth map. Set $\tau _{H}^{I}=f_{kk}^{I}-\zeta ^{k}f_{k}^{I}$.
Then one has
$$
\bigtriangleup _{H}e(f)-f_{i}^{I}\tau _{H,i}^{I}-f_{\alpha
}^{I}\tau _{H,\alpha }^{I}\geq -C_{\varepsilon
}e_{H}(f)-\varepsilon
e_{V}(f)+(f_{ik}^{I})^{2}+\frac{1}{2}(f_{\alpha k}^{I})^{2}
\tag{3.31}
$$
for any given $\varepsilon >0$, where $C_{\varepsilon }$ is a
positive number depending only on $\varepsilon$ and
$$
\sup_{M,i,j,k,\alpha }\{|\zeta ^{k}|,|\zeta _{,i}^{k}|,|\zeta
_{,\alpha }^{k}|,|T_{ij}^{\alpha }|,|T_{ij,k}^{\alpha
}|,|R_{kik}^{j}|,|R_{k\alpha k}^{j}|\}.
$$
In particular, we have
$$
\bigtriangleup _{H}e(f)-f_{i}^{I}\tau _{H,i}^{I}-f_{\alpha
}^{I}\tau _{H,\alpha }^{I}\geq -C_{\varepsilon }e(f).  \tag{3.32}
$$
\endproclaim
\demo{Proof} For any $\varepsilon >0$, we deduce, by Schwarz
inequality, that
$$
\aligned &f_{i}^{I}\zeta
_{,i}^{k}f_{k}^{I}+f_{i}^{I}f_{j}^{I}R_{kik}^{j}\geq
-C_{1}e_{H}(f), \\
&\zeta ^{k}f_{i}^{I}f_{\alpha }^{I}T_{ki}^{\alpha
}+f_{i}^{I}f_{\alpha }^{I}T_{ik,k}^{\alpha }+f_{\alpha }^{I}\zeta
_{,\alpha }^{k}f_{k}^{I}+f_{\alpha }^{I}f_{j}^{I}R_{k\alpha
k}^{j}\geq
-C_{2}(\varepsilon )e_{H}(f)-\varepsilon e_{V}(f_{t}), \\
&2f_{i}^{I}f_{\alpha k}^{I}T_{ik}^{\alpha }\geq
-C_{3}e_{H}(f)-\frac{1}{2} (f_{\alpha k}^{I})^{2},
\endaligned\tag{3.33}
$$
for some positive constants $C_{1}$, $C_{2}(\varepsilon )$ and
$C_{3}$. Since $(N,h)$ has non-positive sectional curvature, we
have
$$
f_{i}^{I}f_{k}^{K}\widetilde{R}_{KJL}^{I}f_{i}^{J}f_{k}^{L}+f_{\alpha
}^{I}f_{k}^{K}\widetilde{R}_{KJL}^{I}f_{\alpha }^{J}f_{k}^{L}\leq
0\text{.} \tag{3.34}
$$
From (3.30), (3.33), (3.34), we obtain (3.31) and thus (3.32)
too.\qed
\enddemo

We will also need similar commutation relations as (3.17) and
(3.21) for maps from the product manifold $M\times (0,\delta)$.
Here the product manifold $M\times (0,\delta)$ is endowed with the
direct sum connection of $\nabla ^{\frak{B}}$ on $M$ and the
trivial connection on $(0,\delta)$. Now let $f:M\times
(0,\delta)\rightarrow N$ be a smooth map. Write
$$
f^{\ast }\widetilde{\omega }=f_{A}^{I}\omega ^{A}+f_{t}^{I}dt.
\tag{3.35}
$$
Taking the exterior derivative of (3.35), one has
$$
Df_{A}^{I}\wedge \omega ^{A}+Df_{t}^{I}\wedge dt+\frac{1}{2}
f_{C}^{I}T_{AB}^{C}\omega ^{A}\wedge \omega ^{B}=0  \tag{3.36}
$$
where
$$
\aligned Df_{A}^{I}&=df_{A}^{I}-f_{B}^{I}\omega
_{A}^{B}+f_{A}^{K}\widetilde{\omega }
_{K}^{I}=f_{AB}^{I}\omega ^{B}+f_{At}^{I}dt \\
Df_{t}^{I}&=df_{t}^{I}+f_{t}^{K}\widetilde{\omega
}_{K}^{I}=f_{tA}^{I}\omega ^{A}+f_{tt}^{I}dt.
\endaligned\tag{3.37}
$$
Consequently $\{f_{AB}^{I}\}$ satisfy (3.17) and
$$
f_{At}^{I}=f_{tA}^{I}.  \tag{3.38}
$$
Similarly taking derivative of the first equation in (3.37) gives
$$
Df_{AB}^{I}\wedge \omega ^{B}+Df_{At}^{I}\wedge
dt=-f_{D}^{I}\Omega _{A}^{D}+f_{A}^{K}\widetilde{\Omega
}_{K}^{I}-f_{AD}^{I}T^{D}  \tag{3.39}
$$
where
$$
\aligned Df_{AB}^{I}&=df_{AB}^{I}-f_{CB}^{I}\omega
_{A}^{C}-f_{AC}^{I}\omega _{B}^{C}+f_{AB}^{K}\widetilde{\omega
}_{K}^{I}=f_{ABC}^{I}\omega
^{C}+f_{ABt}^{I}dt \\
Df_{At}^{I}&=df_{At}^{I}-f_{Bt}^{I}\omega
_{A}^{B}+f_{At}^{K}\widetilde{ \omega }_{K}^{I}=f_{AtC}^{I}\omega
^{C}+f_{Att}^{I}dt.
\endaligned\tag{3.40}
$$
Clearly $\{f_{ABC}^{I}\}$ satisfy (3.21) and
$$
f_{AtB}^{I}-f_{ABt}^{I}=-f_{A}^{K}\widetilde{R}_{KJL}^{I}f_{t}^{J}f_{B}^{L}.
\tag{3.41}
$$
\vskip 0.3 true cm

\heading{\bf 4. Subelliptic harmonic maps and their heat flows}
\endheading
\vskip 0.3 true cm

For a map $f:\left( M^{m+p},H,g_{H};g\right) \rightarrow \left(
N^{n},h\right) $, besides the usual energy $E(f)$, we have the
following two partial energies:
$$
E_{H}(f)=\int_{M}e_{H}(f)dv_{g}=\frac{1}{2}\int_{M}\langle
df(e_{i}),df(e_{i})\rangle dv_{g}  \tag{4.1}
$$
and
$$
E_{V}(f)=\int_{M}e_{V}(f)dv_{g}=\frac{1}{2}\int_{M}\langle
df(e_{\alpha }),df(e_{\alpha })\rangle dv_{g} \tag{4.2}
$$
where the integrands in the second equality of (4.1) (resp. (4.2))
are summed over the range of the index $i$ (resp. $\alpha$). The
partial energies $E_{H}(f)$ and $E_{V}(f)$ are called horizontal
and vertical energies respectively. Clearly
$$
E(f)=E_{H}(f)+E_{V}(f).
$$
\definition{Definition 4.1} A map $f:\left( M,H,g_{H};g\right) \rightarrow
\left( N,h\right) $ is referred to as a subelliptic harmonic map
if it is a critical point of the energy $E_{H}(f)$.
\enddefinition

\proclaim{Proposition 4.1} Let $\{f_{t}\}_{|t|<\varepsilon }$ be a
family of maps from $\left( M,H,g_{H};g\right) \ $to $\left(
N,h\right) $ with $f_{0}=f$ and $\frac{
\partial f_{t}}{\partial t}\mid _{t=0}=\nu \in \Gamma \left( f^{-1}TN\right)
$. Suppose the variation vector field $\nu $ has compact support.
Then
$$
\frac{dE_{H}(f_{t})}{dt}\mid _{t=0}=-\int_{M}\langle \nu ,\tau
_{H}(f)\rangle dv_{g}  \tag{4.3}
$$
where $\tau _{H}(f)=\beta (e_{i},e_{i})-df(\zeta )$ is called the
subelliptic tension field of $f$.
\endproclaim
\demo{Proof} We shall denote by $F:M\times \left( -\varepsilon
,\varepsilon \right) \rightarrow N$ the map defined by
$F(x,t)=f_{t}(x)$. Let $\widetilde{\nabla } ^{F}$ be the pull-back
connection of $\widetilde{\nabla }$ by $F$. Since $
\widetilde{\nabla }$ is torsion-free, we have
$$
\widetilde{\nabla }_{\frac{\partial }{\partial t}}^{F}dF(X)=\widetilde{%
\nabla }_{X}^{F}dF(\frac{\partial }{\partial t})  \tag{4.4}
$$
for any $X\in TM$ (cf. [EL], page 14). Applying (4.1) to $f_{t}$
and using (4.4), we derive that
$$
\aligned \frac{d}{dt}E_{H}(f_{t})|_{t=0}=&\int_{M}\langle
\widetilde{\nabla }_{
\frac{\partial }{\partial t}}^{F}dF(e_{i}),dF(e_{i})\rangle dv_{g}|_{t=0} \\
=&\int_{M}\langle \widetilde{\nabla }_{e_{i}}^{f}\nu
,df(e_{i})\rangle dv_{g}\\
=&\int_{M}{\big (}e_{i}\langle \nu ,df(e_{i})\rangle -\langle \nu
,
\widetilde{\nabla }_{e_{i}}^{f}df(e_{i})\rangle {\big )}dv_{g} \\
=&\int_{M}{\big (}e_{i}\langle \nu ,df(e_{i})\rangle
-\langle \nu ,df(\nabla _{e_{i}}^{\frak{B}}e_{i})\rangle {\big)}dv_{g}\\
&-\int_{M}\langle \nu ,\beta (e_{i},e_{i})\rangle dv_{g},
\endaligned \tag{4.5}
$$
where the terms with the index $i$ are summed over $1\leq i \leq
m$. Set $\theta (X)=\langle \nu ,df\circ \pi _{H}(X)\rangle $ for
any $X\in TM$. The codifferential of $\theta $ is given by
$$
\aligned
\delta \theta &=-(\nabla _{e_{A}}^{R}\theta )(e_{A}) \\
&=-{\big (} e_{A}(\theta (e_{A}))-\theta (\nabla _{e_{A}}^{R}e_{A}){\big )} \\
&=-{\big (}e_{i}\theta (e_{i})-\theta (\pi _{H}(\nabla
_{e_{i}}^{R}e_{i}){\big )} +\theta {\big (}\pi _{H}(\nabla
_{e_{\alpha}}^{R}e_{\alpha }){\big)} \\
&=-{\big (} e_{i}\theta (e_{i})-\theta (\nabla _{e_{i}}^{\frak{B}
}e_{i}){\big )} +\theta (\zeta )
\endaligned\tag{4.6}
$$
where $\zeta =\pi _{H}(\nabla _{e_{\alpha }}^R e_{\alpha })$ (a
sum w.r.t. $\alpha$). It follows from (4.6) and the divergence
theorem that
$$
\int_{M}{\big (}e_{i}\langle \nu ,df(e_{i})\rangle -\langle \nu
,df(\nabla _{e_{i}}^{\frak{B}}e_{i})\rangle {\big )}
dv_{g}=\int_{M}\langle \nu ,df(\zeta )\rangle dv_{g}. \tag{4.7}
$$
By (4.5) and (4.7), we obtain%
$$
\aligned \frac{d}{dt}E_{H}(f_{t})|_{t=0}&=\int_{M}\langle \nu
,df(\zeta )\rangle
dv_{g}-\int_{M}\langle \nu ,\beta (e_{i},e_{i})\rangle dv_{g} \\
&=-\int_{M}\langle \nu ,\tau _{H}(f)\rangle dv_{g}.
\endaligned
$$
\qed
\enddemo

\proclaim{Corollary 4.2} A map $f:(M,H,g_{H};g)\rightarrow (N,h)\
$is a subelliptic harmonic map if
and only if it satisfies the Euler-Lagrange equation%
$$
\tau _{H}(f)=0.  \tag{4.8}
$$
\endproclaim
\remark{Remark 4.1} If $f:\left( M,H,g_{H};g\right) \rightarrow R$
is a smooth function, we find from (3.23) that $\tau
_{H}(f)=\bigtriangleup _{H}f$. Therefore $f$ is a subelliptic
harmonic function if and only if $\bigtriangleup _{H}f=0$.
\endremark

We will introduce a subelliptic heat flow for maps from a
sub-Riemannian manifold $(M,H,$ $g_{H};g)$ to a Riemannian
manifold $(N,h)$ in order to find subelliptic harmonic maps
between these manifolds. Henceforth we assume that both $M$ and
$N$ are compact. As in the theory of harmonic maps, our strategy
to solve (4.8) is to deform a given smooth map $\varphi
:M\rightarrow N$ along the gradient flow of the energy $E_{H}$.
This is equivalent to solving the following subelliptic harmonic
map heat flow:
$$
\cases
\frac{\partial f}{\partial t}=\tau _{H}(f) \\
f|_{t=0}=\varphi
\endcases\tag{4.9}
$$
where $\tau _{H}(f(\cdot ,t))$ is the subelliptic tension field of
$f(\cdot ,t):\left( M,H,g_{H};g\right) \rightarrow \left(
N,h\right) $.

Now we want to give the explicit formulations for both (4.8) and
(4.9), which are convenient for proving the existence theory. In
view of the Nash embedding theorem, one can always assume that
$\frak{I}:(N,h)\rightarrow (R^{K},g_{E})$ is an isometric
embedding in some Euclidean space, where $ g_{E}$ denotes the
standard Euclidean metric. Let $\widetilde{\nabla }$ and $D$
denote the Riemannian connections of $\left( N,h\right) $ and
$\left( R^{K},g_{E}\right) $ respectively. The second fundamental
form of $\frak{I}$ with respect to $(\widetilde{\nabla },D)$ is
$$
\beta (\frak{I};\widetilde{\nabla
},D)(Z,W)=D_{W}d\frak{I}(Z)-d\frak{I}(\widetilde{\nabla }_{W}Z)
\tag{4.10}
$$
where $Z,W$ are any vector fields on $N$. Recall that for a map
$f:(M,\nabla ^{B})\rightarrow (N,\widetilde{\nabla })$, we have
defined its second fundamental form $\beta (f;\nabla
^{\frak{B}},\widetilde{\nabla })$ by (3.12). Applying the
composition formula for second fundamental forms (see Proposition
2.20 on page 16 of [EL]) to the maps $f:\left( M,\nabla
^{\frak{B}} \right) \rightarrow (N,\widetilde{\nabla })$ and
$\frak{I}:(N,\widetilde{\nabla })\rightarrow (R^{K},D)$, we have
$$
\beta (\frak{I}\circ f;\nabla ^{\frak{B}},D)(\cdot ,\cdot
)=d\frak{I}{\big (}\beta (f;\nabla ^{\frak{B}},\widetilde{{\nabla
}} )(\cdot ,\cdot ){\big )}+\beta (\frak{I} ;\widetilde{\nabla
},D){\big (}df(\cdot ),df(\cdot ){\big)}.\tag{4.11}
$$
For simplicity, we shall identify $N$ with $\frak{I}(N)$, and
write $\frak{I}\circ f$ as $u$, which is a map from $M$ to
$R^{K}$. Set
$$
\tau _{H}(u;\nabla ^{\frak{B}},D)=\sum_{i}\beta (u;\nabla
^{\frak{B}} ,D)(e_{i},e_{i})-du(\zeta ).  \tag{4.12}
$$
It follows from (4.11), (4.12) that
$$
\tau _{H}(u;\nabla ^{\frak{B}},D)-tr_{g}\beta (\frak{I};
\widetilde{\nabla },D){\big (}df_{H},df_{H}{\big )}=d\frak{I}{\big
(}\tau _{H}(f){\big)} . \tag{4.13}
$$

By compactness of $N$, there exists a tubular neighborhood $B(N)$
of $N$ in $ R^{K}$ which can be realized as a submersion $\Pi
:B(N)\rightarrow N$ over $ N $. Actually the projection map $\Pi $
is simply given by mapping any point in $B(N)$ to its closest
point in $N$. Clearly its differential $d\Pi
:T_{y}R^{K}\rightarrow T_{y}R^{K}$ when evaluated at a point $y\in
N$ is given by the identity map when restricted to the tangent
space $TN$ of $N$ and maps all the normal vectors to $N$ to the
zero vector. Since $\Pi \circ \frak{I}=\frak{I}:N\hookrightarrow
R^{K}$ and $\beta (\frak{I}; \widetilde{\nabla },D)$ is normal to
$N$, we have
$$
\beta (\frak{I};\widetilde{\nabla },D)(\cdot ,\cdot )=d\Pi {\big
(} \beta (\frak{I};\widetilde{\nabla },D)(\cdot ,\cdot
){\big)}+\beta ( \Pi ;D,D{\big )}(d\frak{I},d\frak{I})
$$
and thus
$$
\beta (\frak{I};\widetilde{\nabla },D)(\cdot ,\cdot )=\beta (\Pi ;
D,D)(d\frak{I},d\frak{I}).  \tag{4.14}
$$
Let $\{y^{a}\}_{1\leq a\leq K}$ be the natural Euclidean
coordinate system of $R^{K}$. Set $u^{a}=y^{a}\circ u$, $\Pi
^{a}=y^{a}\circ \Pi $. From (4.12), Remark 4.1 and (4.14), we have
$$
\tau _{H}(u;\nabla ^{\frak{B}},D)=\bigtriangleup
_{H}u^{a}\frac{\partial }{\partial y^{a}},  \tag{4.15}
$$
and
$$
\aligned tr_{g}\beta (\frak{I};\widetilde{\nabla
},D)(df_{H},df_{H})&=tr_{g}\beta (\Pi
;D,D)(du_{H},du_{H}) \\
&=\Pi _{bc}^{a}\langle \nabla ^{H}u^{b},\nabla ^{H}u^{c}\rangle
\frac{
\partial }{\partial y^{a}}
\endaligned\tag{4.16}
$$
where $\Pi _{bc}^{a}=\frac{\partial ^{2}\Pi ^{a}}{\partial
y^{b}\partial y^{c}}$. Consequently (4.13), (4.15) and (4.16)
imply that
$$
d\frak{I}{\big (}\tau _{H}(f){\big)} ={\big(}\bigtriangleup
_{H}u^{a}-\Pi _{bc}^{a}\langle \nabla ^{H}u^{b},\nabla
^{H}u^{c}\rangle {\big )\frac{
\partial }{\partial y^{a}}.}  \tag{4.17}
$$
Thus $f$ is a subelliptic harmonic map if and only if $u=(u^{a}):M
\rightarrow R^{K}$ satisfies%
$$
\bigtriangleup _{H}u^{a}-\Pi _{bc}^{a}\langle \nabla
^{H}u^{b},\nabla ^{H}u^{c}\rangle =0,\ \ 1\leq a,b,c\leq K.
\tag{4.18}
$$

Inspired by the above explicit formulation for $\tau _{H}(f)$, we
will establish the fact that in order to solve (4.9), it suffices
to solve the following system
$$
\cases \frac{\partial u^{a}}{\partial t}=\bigtriangleup
_{H}u^{a}-\Pi
_{bc}^{a}\langle \nabla ^{H}u^{b},\nabla ^{H}u^{c}\rangle , \\
u^{a}|_{t=0}=\varphi ^{a} \endcases \tag{4.19}
$$
where $\varphi ^{a}=y^{a}\circ \varphi $. Let us define a map
$\rho :B(N)\rightarrow R^{K}$ by
$$
\rho (y)=y-\Pi (y),\quad y\in B(N).
$$
Clearly, $\rho (y)$ is normal to $N$ and $\rho (y)=0$ if and only
if $y\in N$.

\proclaim{Lemma 4.3} Let $u(x,t)$ $=(u^{a}(x,t))$ $((x,t)\in
M\times \lbrack 0,\delta))$ be a solution of $(4.19)$ with initial
condition $\varphi $ $=(\varphi ^{a}):M\rightarrow R^{K}$. Then
the quantity
$$
\int_{M}|\rho (u(x,t))|^{2}dv_{g}
$$
is a nonincreasing function of $t$. In particular, if $\varphi
(M)\subset N$, then $u(x,t)\in N$ for all $(x,t)\in M\times
\lbrack 0,\delta)$.
\endproclaim
\demo{Proof} Since $\rho (y)=y-\Pi (y)$, we have
$$
\rho _{b}^{a}=\delta _{b}^{a}-\Pi _{b}^{a}  \tag{4.20}
$$
and
$$
\rho _{bc}^{a}=-\Pi _{bc}^{a}  \tag{4.21}
$$
where $\rho _{b}^{a}=\frac{\partial \rho ^{a}}{\partial y^{b}}$
and $\rho _{bc}^{a}=\frac{\partial ^{2}\rho ^{a}}{\partial
y^{b}\partial y^{c}}$. By applying the composition law ([EL]) to
the maps $u_{t}:(M,\nabla ^{\frak{B}})\rightarrow (B(N),D)$ and
$\rho :(B(N),D)\rightarrow (R^{K},D)$, we have
$$
\bigtriangleup _{H}\rho (u)=d\rho (\bigtriangleup
_{H}u)+tr_{g}\beta (\rho ;D,D)(du_{H},du_{H}).  \tag{4.22}
$$
It follows from (4.20), (4.21), (4.22) and (4.19) that
$$
\aligned (\bigtriangleup _{H}\rho (u))^{a}&=\rho
_{b}^{a}\bigtriangleup _{H}u^{b}+\rho
_{bc}^{a}\langle \nabla ^{H}u^{b},\nabla ^{H}u^{c}\rangle \\
&=\bigtriangleup _{H}u^{a}-\Pi _{b}^{a}\bigtriangleup
_{H}u^{b}-\Pi
_{bc}^{a}\langle \nabla ^{H}u^{b},\nabla ^{H}u^{c}\rangle \\
&=\frac{\partial u^{a}}{\partial t}-\Pi _{b}^{a}\bigtriangleup _{H}u^{b} \\
&=\rho _{b}^{a}\frac{\partial u^{b}}{\partial t}+\Pi
_{b}^{a}(\frac{\partial u^{b}}{\partial t}-\bigtriangleup
_{H}u^{b}).
\endaligned\tag{4.23}
$$
Since $d\Pi (\frac{\partial u}{\partial t}-\bigtriangleup _{H}u)$
is tangent to $N$ and $\rho (u)$ is normal to $N$, we find from
(4.23) that
$$
\rho ^{a}(u)(\bigtriangleup _{H}\rho (u))^{a}=\rho ^{a}(u)\rho
_{b}^{a}(u)\frac{\partial u^{b}}{\partial t}.  \tag{4.24}
$$
Using (4.24), (2.5), we deduce that
$$
\aligned \frac{\partial }{\partial t}\int_{M}(\rho
^{a}(u))^{2}dv_{g}&=2\int_{M}\rho
^{a}\rho _{b}^{a}(u)\frac{\partial u^{b}}{\partial t}dv_{g} \\
&=2\int_{M}\rho ^{a}(u)(\bigtriangleup _{H}\rho (u))^{a}dv_{g} \\
&=-2\int_{M}|\nabla _{H}\rho (u)|^{2}dv_{g} \\
&\leq 0
\endaligned
$$
which proves this lemma.\qed
\enddemo

In terms of (4.17) and Lemma 4.3, we conclude that
\proclaim{Theorem 4.4} Let $\varphi :M\rightarrow N\subset R^{K}$
be a smooth map given by $\varphi =(\varphi ^{1},...,\varphi
^{K})$ in the Euclidean coordinates. If $ u:M\times \lbrack
0,\delta)\rightarrow N$ is a solution of the following system
$$
\frac{\partial u^{a}}{\partial t}=\bigtriangleup _{H}u^{a}-\Pi
_{bc}^{a}\langle \nabla ^{H}u^{b},\nabla ^{H}u^{c}\rangle ,\quad
1\leq a\leq K,
$$
with initial condition $(u^{a}(x,0))=(\varphi ^{a}(x))$ for all
$x\in M$,
then $u$ solves the subelliptic heat flow%
$$
\frac{\partial u}{\partial t}=\tau _{H}(u)
$$
with initial condition $u(x,0)=\varphi (x)$.
\endproclaim

A general version of the second variation formula for $E_{H}$ is
useful for our purpose. Although its derivation is routine, we now
derive this formula for the convenience of the readers.

\proclaim{Proposition 4.5} Let $F:\left( M,H,g_{H};g\right) \times
(-\varepsilon ,\varepsilon )\rightarrow N$ be a family of maps
with $F(\cdot ,0)=f$ and $\frac{\partial F}{\partial t}\mid
_{t=0}=\nu \in \Gamma \left( f^{-1}TN\right) $. Then
$$
\aligned \frac{d^{2}E_{H}{\big (}F(\cdot ,t){\big
)}}{dt^{2}}|_{t=0}
=&-\int_{M}\langle \xi ,\tau _{H}(f)\rangle dv_{g} \\
&+\int_{M}{\big\{}\langle \widetilde{\nabla
}_{e_{i}}\nu,\widetilde{\nabla }_{e_{i}}\nu\rangle
-\widetilde{R}{\big(}df(e_{i}),\nu,df(e_{i}),\nu{\big)}{\big\}}dv_{g}\endaligned
$$
where $\xi =\widetilde{\nabla }_{\frac{\partial }{\partial
t}}dF(\frac{\partial }{\partial t})|_{t=0}\in \Gamma (f^{-1}TN)$.
\endproclaim
\demo{Proof} At each $t$, we compute
$$
\aligned \frac{\partial e_{H}{\big (}F(\cdot ,t){\big )}}{\partial
t}&=\langle \widetilde{\nabla }_{\frac{\partial }{\partial
t}}^{F}dF(e_{i}),dF(e_{i})
\rangle \\
&=\langle \widetilde{\nabla }_{e_{i}}^{F}dF(\frac{\partial
}{\partial t} ),dF(e_{i})\rangle
\endaligned \tag{4.25}
$$
and
$$
\frac{\partial ^{2}e_{H}}{\partial t^{2}}=\langle
\widetilde{\nabla }_{\frac{
\partial }{\partial t}}^{F}\widetilde{\nabla }_{e_{i}}^{F}dF(\frac{\partial
}{\partial t}),dF(e_{i})\rangle +\langle \widetilde{\nabla
}_{e_{i}}^{F}dF( \frac{\partial }{\partial t}),\widetilde{\nabla
}_{e_{i}}^{F}dF(\frac{\partial }{\partial t})\rangle  \tag{4.26}
$$
where $\{e_{i}\}$ is a local orthonormal frame field for
$(H,g_{H})$. Note that
$$
\aligned \widetilde{\nabla }_{\frac{\partial }{\partial
t}}^{F}\widetilde{\nabla } _{X}^{F}&=\widetilde{\nabla
}_{X}^{F}\widetilde{\nabla }_{\frac{\partial }{
\partial t}}^{F}+\widetilde{R}^{F}(\frac{\partial }{\partial t},X)+
\widetilde{\nabla }_{[\frac{\partial }{\partial t},X]}^{F} \\
&=\widetilde{\nabla }_{X}^{F}\widetilde{\nabla }_{\frac{\partial
}{\partial t} }^{F}+\widetilde{R}(dF(\frac{\partial }{\partial
t}),dF(X)).
\endaligned\tag{4.27}
$$
From (4.26) and (4.27), we obtain%
$$
\aligned \frac{\partial ^{2}e_{H}}{\partial t^{2}}=&\langle
\widetilde{\nabla }_{e_{i}} \widetilde{\nabla }_{\frac{\partial
}{\partial t}}dF(\frac{\partial }{
\partial t}),dF(e_{i})\rangle +\langle \widetilde{R}{\big(}dF(\frac{\partial }{
\partial t}),dF(e_{i}){\big)}dF(\frac{\partial }{\partial t}),dF(e_{i})\rangle \\
&+\langle \widetilde{\nabla }_{e_{i}}dF(\frac{\partial }{\partial
t}), \widetilde{\nabla }_{e_{i}}dF(\frac{\partial }{\partial
t})\rangle
\endaligned
$$
and thus
$$
\aligned &\frac{d^{2}E_{H}(F)}{\partial
t^{2}}|_{t=0}\\
&=\int_{M}\langle \widetilde{ \nabla
}_{e_{i}}^{F}\widetilde{\nabla }_{\frac{\partial }{\partial
t}}^{F}dF(
\frac{\partial }{\partial t}),dF(e_{i})\rangle dv_{g}\mid _{t=0} \\
&\phantom{ab}+\int_{M}\{\langle \widetilde{\nabla
}_{e_{i}}dF(\frac{\partial }{\partial t} ),\widetilde{\nabla
}_{e_{i}}dF(\frac{\partial }{\partial t})\rangle +\langle
\widetilde{R}(dF(\frac{\partial }{\partial t}),dF(e_{i}))dF(\frac{
\partial }{\partial t}),dF(e_{i})\rangle \}dv_{g}\mid _{t=0} \\
&=\int_{M}\langle \widetilde{\nabla }_{e_{i}}\xi ,df(e_{i})\rangle
dv_{g}+\int_{M}\{\langle \widetilde{\nabla
}_{e_{i}}\nu,\widetilde{\nabla } _{e_{i}}\nu\rangle
-\widetilde{R}(df(e_{i}),\nu,df(e_{i}),\nu)\}dv_{g}
\endaligned\tag{4.28}
$$
where $\xi =\widetilde{\nabla }_{\frac{\partial }{\partial
t}}dF(\frac{\partial }{\partial t})|_{t=0}\in \Gamma (f^{-1}TN)$.
By (4.5) and (4.7), we have
$$
\int_{M}\langle \widetilde{\nabla }_{e_{i}}^{f}\xi
,df(e_{i})\rangle dv_{g}=-\int_{M}\langle \xi ,\tau _{H}(f)\rangle
dv_{g}.  \tag{4.29}
$$
In terms of (4.28) and (4.29), we complete the proof of this
proposition.\qed
\enddemo

\proclaim{Corollary 4.6} Suppose $f:M\times \lbrack
0,\delta)\rightarrow N$ is a solution of the subelliptic harmonic
map heat flow $\partial f/\partial t=\tau _{H}{\big (} f(\cdot
,t){\big )}$ for $t\in \lbrack 0,\delta)$. Then
$$
\frac{d^{2}E_{H}{\big (}f(\cdot ,t){\big
)}}{dt^{2}}=2\int_{M}{\big \{} \langle \widetilde{\nabla
}_{e_{i}}\tau _{H}(f),\widetilde{\nabla }_{e_{i}}\tau
_{H}(f)\rangle -\widetilde{R}(df(e_{i}),\tau
_{H}(f_{t}),df(e_{i}),\tau _{H}(f)){\big \}dv}_{g}.
$$
\endproclaim
\demo{Proof} Applying Proposition 4.5 to $\{f(\cdot ,t)\}$ at each
$t\in \lbrack 0,\delta)$, we get
$$
\aligned \frac{d^{2}E_{H}{\big(}f(\cdot
,t){\big)}}{dt^{2}}=&-\int_{M}\langle \widetilde{\nabla }_{
\frac{\partial }{\partial t}}\tau _{H}(f),\tau _{H}(f)\rangle dv_{g} \\
&+\int_{M}{\big\{}\langle \widetilde{\nabla }_{e_{i}}\tau
_{H}(f),\widetilde{\nabla }_{e_{i}}\tau _{H}(f)\rangle
-\widetilde{R}{\big(}df(e_{i}),\tau _{H}(f),df(e_{i}),\tau
_{H}(f){\big)}{\big\}}dv_{g}. \endaligned\tag{4.30}
$$
Note that Proposition 4.1 gives
$$
\frac{dE_{H}{\big (}f(\cdot ,t){\big )}}{dt}=\int_{M}\langle
\frac{
\partial f}{\partial t},\tau _{H}(f)\rangle dv_{g}=-\int_{M}|\tau
_{H}(f)|^{2}dv_{g}.
$$
Consequently
$$
\frac{d^{2}E_{H}{\big (}f(\cdot ,t){\big
)}}{dt^{2}}=-2\int_{M}\langle \widetilde{\nabla }_{\frac{\partial
}{\partial t}}\tau _{H}(f),\tau _{H}(f)\rangle dv_{g}.  \tag{4.31}
$$
This corollary follows immediately from (4.30) and (4.31).\qed
\enddemo

\heading{\bf 5. Existence of Subelliptic Harmonic Maps}
\endheading
\vskip 0.3 true cm

{\bf 5.1 Short-time Existence} \vskip 0.3 true cm

For bounded functions $f:M\times \lbrack 0,\delta)\rightarrow R$
and $\psi :M\rightarrow R$, let us consider the subelliptic heat
flow
$$
\cases(\bigtriangleup _{H}-\frac{\partial }{\partial t})w=f(x,t), \\
w|_{t=0}=\psi .
\endcases\tag{5.1}
$$
By Duhamel's principle, we know that one solution of (5.1) is given by%
$$
w(x,t)=\int_{M}K(x,y,t)\psi
(y)dv_{g}-\int_{0}^{t}\int_{M}K(x,y,t-s)f(y,s)dv_{g}(y)ds.
\tag{5.2}
$$

First we establish the following short-time existence theorem.
\proclaim{Theorem 5.1}
Let $(M^{m+d},H,g_{H},g)$ be a compact sub-Riemannian manifold, and $%
(N^{n},h)\subset R^{K}$ be a compact submanifold with the induced
Euclidean metric. For any smooth map $\varphi :M\rightarrow N$,
there exists $\delta_{0}>0$ such that the subelliptic harmonic map
heat flow with initial condition
$$
\cases (\bigtriangleup _{H}-\frac{\partial }{\partial
t})u^{a}(x,t)=\Pi
_{bc}^{a}\langle \nabla ^{H}u^{b},\nabla ^{H}u^{c}\rangle \\
\\
u^{a}(x,0)=\varphi ^{a}(x) ,\quad 1\leq a,b,c\leq K.
\endcases
$$
admits a smooth solution on $M\times \lbrack 0,\delta_{0})$, where
$\delta_{0}\ $is a constant depending only on $\sup_{M}e(\varphi
)$ and geometric quantities of both $M$ and $N$.
\endproclaim
\demo{Proof} Writing $u=(u^{a}(x,t))_{1\leq a\leq K}$, the
subelliptic harmonic map heat flow may be expressed as
$$
\cases
(\bigtriangleup _{H}-\frac{\partial }{\partial t})u=F(x,t) \\
u(x,0)=\varphi (x)
\endcases  \tag{5.3}
$$
where $F(x,t)=(\Pi _{bc}^{a}\langle \nabla ^{H}u^{b},\nabla
^{H}u^{c}\rangle )$ depends on the unknown solution $u$ itself. In
terms of (5.2), we can define a sequence of approximate solutions
for (5.3) inductively as follows:
$$
\aligned
u_{0}(x,t)&=\int_{M}K(x,y,t)\varphi (y)dv_{g}(y) \\
u_{k}(x,t)&=u_{0}(x,t)-\int_{0}^{t}\int_{M}K(x,y,t-s)F_{k-1}(y,s)dv_{g}(y)ds
\endaligned\tag{5.4}
$$
where
$$
F_{k-1}(y,s)={\big (}\Pi _{bc}^{a}\langle \nabla
^{H}u_{k-1}^{b},\nabla ^{H}u_{k-1}^{c}\rangle (y,s){\big )},\quad
k=1,2,3,\cdots .  \tag{5.5}
$$
Clearly $u_{0}$ and $u_{k}:M\rightarrow R^{K}$ satisfy
respectively
$$
\cases
\left( \triangle _{H}-\frac{\partial }{\partial t}\right) u_{0}=0, \\
u_{0}(x,0)=\varphi (x)%
\endcases \tag{5.6}
$$
and
$$
\cases
\left( \triangle _{H}-\frac{\partial }{\partial t}\right)
u_{k}=F_{k-1}(x,t)
\\
u_{k}(x,0)=\varphi (x) ,\quad k=1,2,\cdots .\endcases \tag{5.7}
$$
We set
$$
\Lambda =\sup_{B(N),\,a,b,c,d}\{\mid \Pi _{bc}^{a}\mid,\phantom{}
\mid \frac{\partial \Pi _{bc}^{a}}{\partial y^{d}}\mid \}
\tag{5.8}
$$
where $(y^{1},...,y^{K})$ are coordinates of $R^{K}$, and $B(N)$
is the tubular neighborhood of $N$ on which $\Pi $ is defined. Let
us also introduce
$$
p_{k-1}(t)=\sup_{M\times \lbrack 0,t]}\sqrt{e_{H}(u_{k-1})},\quad
k=1,2,\cdots .  \tag{5.9}
$$
which is obviously non-decreasing in $t$. From (5.5) and (5.9), we
have
$$
\sup_{M\times \lbrack 0,t]}\mid F_{k-1}(x,s)\mid \leq \Lambda
p_{k-1}^{2}(t). \tag{5.10}
$$
Note that
$$
\mid u_{0}\mid \leq \Vert \varphi \Vert _{C^{0}}=\sup_{x\in
M}\sqrt{\sum_{a=1}^{K}(\varphi ^{a}(x))^{2}},  \tag{5.11}
$$
since $\int_{M}K(x,y,t)dy=1$. Here and afterwards, $\Vert \cdot
\Vert _{C^{0}}$ denotes the $C^{0}$-norm of functions or tensor
fields on $M$. From (5.4), (5.10) and (5.11), we derive that
$$
\mid u_{k}-u_{0}\mid \leq \Lambda tp_{k-1}^{2}  \tag{5.12}
$$
and
$$
\mid u_{k}\mid \leq \Lambda tp_{k-1}^{2}(t)+\Vert \varphi \Vert
_{C^{0}}. \tag{5.13}
$$
Note that $\tau _{H}(u_{0})=\bigtriangleup _{H}u_{0}$ for the map
$ u_{0}:M\rightarrow R^{K}$. In view of (3.32), we have
$$
\left( \bigtriangleup _{H}-\frac{\partial }{\partial t}\right)
e(u_{0})\geq -Ce(u_{0}),
$$
or equivalently,
$$
\left( \bigtriangleup _{H}-\frac{\partial
}{\partial t}\right) \left( e^{-Ct}e(u_{0})\right) \geq 0.
\tag{5.14}
$$
Consequently the Maximum principle\ (see Lemma 2.4) implies that
$$
e^{-Ct}e(u_{0})\leq e(\varphi)
$$
and thus
$$
p_{0}(t)\leq \sqrt{e^{Ct}e(\varphi )}.  \tag{5.15}
$$
Using Lemma 2.3, (5.4) and (5.10), we may deduce
$$
\aligned |\nabla _{x}^{H}u_{k}(x,t)|&\leq |\nabla
_{x}^{H}u_{0}|+\int_{0}^{t}\int_{M}|\nabla
_{x}^{H}K(x,y,t-s)|\,|F_{k-1}(y,s)|dv_{g}(y) \\
&\leq |\nabla _{x}^{H}u_{0}|+C_1\Lambda p_{k-1}^{2}(t)t^{\beta }
\endaligned
$$
hence implying
$$
p_{k}(t)\leq C_{1}\Lambda p_{k-1}^{2}(t)t^{\beta }+p_{0}(t).
\tag{5.16}
$$
For any $0<\varepsilon <1$, by choosing $\delta$ sufficiently
small, (5.15) yields that
$$
C_{1}\Lambda \delta^{\beta }p_{0}(\delta)\leq C_{1}\Lambda
\delta^{\beta }\sqrt{ e^{C\delta}e(\varphi )}\leq
\frac{\varepsilon }{4}.
$$
By an inductive argument, we get
$$
C_{1}\Lambda \delta^{\beta }p_{k}(\delta)\leq \frac{\varepsilon
}{2} \tag{5.17}
$$
since (5.16) gives
$$
\aligned C_{1}\Lambda \delta^{\beta }p_{k}(\delta) &\leq \left(
C_{1}\Lambda \delta^{\beta
}p_{k-1}(\delta)\right) ^{2}+C_{1}\Lambda \delta^{\beta }p_{0}(\delta) \\
&\leq \frac{\varepsilon }{4}+\frac{\varepsilon }{4} \\
&=\frac{\varepsilon }{2}.
\endaligned
$$
Consequently
$$
p_{k}(\delta)\leq C_{2}\varepsilon \delta^{-\beta }.  \tag{5.18}
$$
We define the following space of functions,
$$
C_{H}^{1}(M,R^{K})=\{f:M\rightarrow R^{K}\mid f\in C^{0},\quad
\nabla ^{H}f\in C^{0}\}
$$
which is endowed with the norm
$$
\Vert f\Vert _{C_{H}^{1}}=\Vert f\Vert _{C^{0}}+\Vert \nabla
^{H}f\Vert _{C^{0}}.
$$
It is known that $(C_{H}^{1}(M,R^{K}),\Vert \cdot \Vert
_{C_{H}^{1}})$ is a Banach space. From (5.13) and (5.18), one has
$$
||u_{k}||_{C_{H}^{1}(M,R^{K})}\leq C_{3}(C_{2},\varepsilon
,\delta).
$$
In terms of (5.12) on $M\times \lbrack 0,\delta)$ and using
(5.17), we deduce that
$$
\aligned
|u_{k}(x,t)-u_{0}(x,t)|&\leq \Lambda \delta p_{k-1}^{2}(\delta) \\
&\leq \frac{\varepsilon ^{2}\delta^{1-2\beta }}{4C_{1}\Lambda }.
\endaligned\tag{5.19}
$$
The validity of the inequality (5.17) depends on choosing a
sufficiently\ small $\delta$. Note also that $1-2\beta >0$. From
(5.19), we find that all maps $u_{k}$ ($k=1,2,...$) will map $M$
into $B(N)$ by choosing both $\varepsilon $ and $\delta$
sufficiently small since $|u_{0}(x,t)-\varphi (x)|$ can be chosen
to be sufficiently small for small $t$ by continuity of $u_{0}$.

Now we want to show that $\{u_{k}(x,t)\}$ form a Cauchy sequence
in $ C_{H}^{1}(M,R^{K})$ for sufficiently small $t$. Let us define
$$
X_{k}(t)=\sup_{M\times \lbrack 0,t]}\{\mid
u_{k}(x,s)-u_{k-1}(x,s)\mid +\mid \nabla _{x}^{H}u_{k}(x,s)-\nabla
_{x}^{H}u_{k-1}(x,s)\mid \} \tag{5.20}
$$
which is a non-decreasing function of $t$. Note that
$$
\aligned F_{k}(x,t)-F_{k-1}(x,t)=&{\big (}\Pi
_{bc}^{a}(u_{k})\langle \nabla ^{H}u_{k}^{b},\nabla
^{H}u_{k}^{c}\rangle -\Pi _{bc}^{a}(u_{k-1})\langle
\nabla ^{H}u_{k-1}^{b},\nabla ^{H}u_{k-1}^{c}\rangle {\big )} \\
=&{\big (}\Pi _{bc}^{a}(u_{k})-\Pi _{bc}^{a}(u_{k-1}){\big
)}\langle \nabla ^{H}u_{k}^{b},\nabla ^{H}u_{k}^{c}\rangle\\
&+{\big (}\Pi _{bc}^{a}(u_{k-1})(\langle \nabla
^{H}u_{k}^{b},\nabla ^{H}u_{k}^{c}\rangle -\langle \nabla
^{H}u_{k-1}^{b},\nabla ^{H}u_{k-1}^{c}\rangle ){\big )}\\
=&{\Big(}{\big (}\Pi _{bc}^{a}(u_{k})-\Pi _{bc}^{a}(u_{k-1}){\big
)} \langle \nabla ^{H}u_{k}^{b},\nabla ^{H}u_{k}^{c}\rangle {\Big
)}\\
&+{\Big(}\Pi _{bc}^{a}(u_{k-1})(\langle \nabla
^{H}u_{k}^{b}-\nabla
^{H}u_{k-1}^{b},\nabla ^{H}u_{k}^{c}\rangle ){\Big )} \\
&+{\Big (}\Pi _{bc}^{a}(u_{k-1})(\langle \nabla
^{H}u_{k-1}^{b},\nabla ^{H}u_{k}^{c}-\nabla ^{H}u_{k-1}^{c}\rangle
{\Big)}.
\endaligned \tag{5.21}
$$
Using (5.18) and the estimate
$$
\mid \Pi _{bc}^{a}(u_{k})-\Pi _{bc}^{a}(u_{k-1})\mid \leq
\Lambda \mid u_{k}-u_{k-1}\mid ,
$$
we may derive from (5.21) that
$$
\aligned \sup_{M\times \lbrack 0,t]}\mid
F_{k}(x,t)-F_{k-1}(x,t)\mid &\leq
C_{4}X_{k}(t){\big (}p_{k}^{2}(t)+p_{k}(t)+p_{k-1}(t){\big )} \\
&\leq C_{5}X_{k}(t)%
\endaligned\tag{5.22}
$$
for any $t\leq \delta$. Consequently we get the following two
estimates
$$
\aligned \mid u_{k}-u_{k-1}\mid &\leq
\int_{0}^{t}\int_{M}K(x,y,t-s)\mid
F_{k-1}(y,s)-F_{k-2}(y,s)\mid dv_{g}(s)ds \\
&\leq C_{5}tX_{k-1}(t)
\endaligned
$$
and
$$
\aligned
&\mid \nabla _{x}^{H}u_{k}-\nabla _{x}^{H}u_{k-1}\mid \\
&\leq \int_{0}^{t}\int_{M}\mid \nabla _{x}^{H}K(x,y,t-s)\mid \cdot
\mid
F_{k-1}(y,s)-F_{k-2}(y,s)\mid dv_{g}(s)ds \\
&\leq C_{6}t^{\beta }X_{k-1}(t)
\endaligned
$$
which imply
$$
X_{k}(t)\leq C_{7}t^{\beta }X_{k-1}(t)  \tag{5.23}
$$
for $k\geq 2$. For $k=1$, using $t<1$, we have from (5.4) and
(5.15) that
$$
\aligned \mid u_{1}(x,t)-u_{0}(x,t)\mid &\leq
\int_{0}^{t}\int_{M}K(x,y,t-s)\mid
F_{0}(y,s)\mid dv_{g}(y)ds \\
&\leq t\Lambda p_{0}^{2}(t) \\
&\leq t\Lambda e^{C}e(\varphi )
\endaligned  \tag{5.24}
$$
and
$$
\aligned \mid \nabla _{x}^{H}u_{1}(x,t)-\nabla
_{x}^{H}u_{0}(x,t)\mid &\leq \int_{0}^{t}\int_{M}\mid \nabla
_{x}^{H}K(x,y,t-s)\mid \mid F_{0}(y,s)\mid
dv_{g}(y)ds \\
&\leq C_{1}t^{\beta }\Lambda p_{0}^{2}(t) \\
&\leq C_{1}t^{\beta }\Lambda e^{C}e(\varphi ).
\endaligned\tag{5.25}
$$
It follows that
$$
X_{1}(t)\leq C_{8}(C_{7}t^{\beta })e(\varphi ).  \tag{5.26}
$$
By iterating (5.23) and using (5.26), we get
$$
X_{k}(t)\leq C_{8}(C_{7}t^{\beta })^{k}e(\varphi ).  \tag{5.27}
$$
We may choose a sufficiently small positive number $\delta_{0}$
such that $\delta_{0}\leq \delta$ and $C_{7}\delta_{0}^{\beta
}<1$. Hence (5.27) implies that for any $ i<j$
$$
\aligned &{\sup }_{\lbrack 0,\delta_{0}]}\parallel u_{i}(\cdot
,t)-u_{j}(\cdot
,t)\parallel _{C_{H}^{1}(M)}\\
&\leq \sum_{k=i+1}^{j}X_{k}(\delta_{0}) \\
&\leq C_{9}\sum_{k=i+1}^{j}(C_{7}\delta_{0}^{\beta })^{k}
\endaligned
$$
which tends to $0$ as $i,j\rightarrow \infty $. Hence there exists
$u\in C^{0}(M\times \lbrack 0,\delta_{0}],B(N))$ with $u(\cdot
,t)\in C_{H}^{1}(M,B(N))$ for each $t\in \lbrack 0,\delta_{0}]$,
such that $u_{k}\rightarrow u$ and $\nabla ^{H}u_{k}\rightarrow
\nabla ^{H}u$ uniformly on $M\times \lbrack 0,\delta_{0}]$.
Consequently
$$
F_{k}(x,t)\rightarrow F(x,t)={\big (}\Pi _{bc}^{a}\left( u\right)
\langle \nabla ^{H}u^{b},\nabla ^{H}u^{c}\rangle {\big )}
$$
and thus (5.4) implies that $u$ is given by%
$$
u(x,t)=\int_{M}K(x,y,t)\varphi
(y)dv_{g}(y)-\int_{0}^{t}\int_{M}K(x,y,t-s)F(y,s)dv_{g}(y)ds.
$$
Clearly $u$ is a weak solution of the subelliptic harmonic map
heat flow. In terms of Theorem 2.1 and Remark 2.1, by a
bootstrapping argument, we find that $u\in C^{\infty }(M\times
(0,\delta_{0}),N)$ satisfies (4.19).\qed
\enddemo

Next we give the following uniqueness theorem.

\proclaim{Theorem 5.2} Let $u$ and $v$ be solutions on $M\times
\lbrack 0,\delta)$ to the subelliptic harmonic map heat flow with
the same initial condition $:$ $u(x,0)=v(x,0)=\varphi (x)$. Then
$u$ and $v$ are identical.
\endproclaim
\demo{Proof} Set $\Psi =\sum_{a=1}^{K}\left( u^{a}-v^{a}\right)
^{2}$. A direct computation gives
$$
\aligned \left( \triangle _{H}-\frac{\partial }{\partial t}\right)
\Psi &=2\sum_{a}(u^{a}-v^{a})\left( \triangle _{H}-\frac{\partial
}{\partial t}\right) (u^{a}-v^{a})+2\sum_{a}\mid \nabla ^{H}(u^{a}-v^{a})\mid ^{2}\\
&=2\sum_{A}(u^{a}-v^{a})(F^{a}(u)-F^{a}(v))+2\sum_{a}\mid \nabla
^{H}(u^{a}-v^{a})\mid ^{2}.
\endaligned\tag{5.28}
$$
For any $0<\delta_{1}<\delta$, we set
$$
p_{\delta_{1}}=\sup_{M\times \lbrack
0,\delta_{1}]}\sqrt{e_{H}(u)},\quad q_{\delta_{1}}=\sup_{M\times
\lbrack 0,\delta_{1}]}\sqrt{e_{H}(v)}.
$$
Writing $F^{a}(u)-F^{a}(v)$ in a similar way as (5.21), one may
get
$$
\mid F^{a}(u)-F^{a}(v)\mid \leq C(\delta_{1},\Lambda
,p_{\delta_{1}},q_{\delta_{1}}){\big (}\Psi
^{\frac{1}{2}}+\sum_{b}\mid \nabla ^{H}(u^{b}-v^{b})\mid {\big )}
\tag{5.29}
$$
on any $[0,\delta_{1}]$ with $\delta_{1}<\delta$, where
$C(\delta_{1},\Lambda ,p_{\delta_{1}},q_{\delta_{1}})$ is a
constant depending on $\delta_{1}$, $\Lambda $, $p_{\delta_{1}}$
and $q_{\delta_{1}}$. It follows immediately from (5.28) and
(5.29) that
$$
\left( \triangle _{H}-\frac{\partial }{\partial t}\right) \Psi
\geq -\widetilde{C}(\delta_{1},\Lambda
,p_{\delta_{1}},q_{\delta_{1}})\Psi
$$
on $[0,\delta_{1}]$ for some positive constant
$\widetilde{C}(\delta_{1},\Lambda
,p_{\delta_{1}},q_{\delta_{1}})$. This implies
$$
\left( \triangle _{H}-\frac{\partial }{\partial t}\right)
(e^{-\widetilde{C} t}\Psi )\geq 0
$$
on $[0,\delta_{1}]$ and thus the maximum principle asserts that
$\Psi =0$ on $ [0,\delta_{1}]$. Since $\delta_{1}$ is arbitrary,
we conclude that $\Psi =0$ on $[0,\delta). $\qed
\enddemo
\vskip 0.5 true cm

{\bf 5.2 Long-time Existence} \vskip 0.4 true cm

We first give a criteria for the long-time existence of the
subelliptic harmonic map heat flow.

\proclaim{Lemma 5.3} Suppose $u=\frak{I}\circ f$ is a solution of
the subelliptic harmonic heat flow on $M\times \lbrack
0,\delta_{\max })$, where $\delta_{\max }$ is the maximal
existence time for the solution $u$. If $\delta_{\max }<\infty $,
then
$$
{\lim\inf}_{t\rightarrow \delta_{\max }}\{\sup_{M}e{\big (}u(\cdot
,t) {\big )}\}=+\infty \text{.}
$$
In other words, if
$$
{\lim\inf}_{t\rightarrow \delta-0}\{\sup_{M}e{\big (}u(\cdot ,t)
{\big )}\}<+\infty
$$
on any $M\times \lbrack 0,\delta)$, where the solution $u$ exists,
then $\delta_{\max }=\infty $ (long-time existence).
\endproclaim
\demo{Proof} Suppose $u$ is a solution of the subelliptic harmonic
map heat flow on $M\times \lbrack 0,\delta_{\max })$ with
$\delta_{\max }<\infty $. We want to prove that
$$
{\lim \inf }_{t\rightarrow \delta_{\max }}\{\sup_{M}e{\big
(}u(\cdot ,t){\big )}\}=+\infty .
$$
Otherwise, there is a sequence $t_{k}\rightarrow \delta_{\max }$
such that $\sup_{M}e{\big (}u(\cdot ,t_{k}){\big )}\leq C_{0}$ for
some positive number $C_{0}$. By Theorem 5.1, there exists a
positive number $\delta(C_{0},M,N)$ depending only on $C_{0}$ and
the geometric quantities of $M$ and $N$ such that the subelliptic
harmonic heat flow admits a solution with $u_{t_{k}}$ as its
initial condition on $[t_{k},t_{k}+\delta(C_{0},M,N))$. Taking a
sufficiently large $k$, the uniqueness in Theorem 5.2 enables us
to obtain a solution on $M\times \lbrack 0,\delta_{\max }+\delta
^{\prime })$ for some positive number $\delta ^{\prime }$. This
contradicts to the assumption that $\delta_{\max }$ is the maximal
existence time.\qed
\enddemo

From now on, we assume that $(N,h)$ has non-positive sectional
curvature. Let $f:M\rightarrow N\ $ be a solution of the
subelliptic harmonic map heat flow on $[0,\delta)$. By (3.32) we
get
$$
(\bigtriangleup _{H}-\frac{\partial }{\partial t})e(f)\geq -Ce(f)
$$
for some constant, that is,
$$
(\bigtriangleup _{H}-\frac{\partial }{\partial t}){\big
(}e^{-Ct}e(f) {\big)}\geq 0.  \tag{5.30}
$$

\proclaim{Lemma 5.4} Let $f:M\rightarrow N\ $ be a solution of the
subelliptic harmonic map heat flow on $[0,\delta)$. Suppose
$(N,h)$ has non-positive sectional curvature. Set $ \alpha =\min
\{R_{0},\sqrt{\delta}\}$, where $R_0$ is given by Lemma 2.4. Then
$$
e{\big (}f(\cdot ,t){\big )}\leq C(\varepsilon _{0})E{\big
(}f(\cdot ,t-\varepsilon _{0}){\big )}
$$
for $t\in \lbrack \varepsilon _{0},\delta)$, where $\varepsilon
_{0}$ is a fixed number in $(0,\frac{\alpha^{2}}{2})$.
\endproclaim
\demo{Proof} Using (5.30) and applying the mean value inequality
in Lemma 2.4 to $ e^{-C(s+t)}e{\big (}(x,s+t){\big )}$ for $t\in $
$(0,\alpha^{2})$ and $ s+t<\delta$, we obtain
$$
e^{-C(s+t)}e(f(x,s+t))\leq Bt^{-\frac{Q}{2}}\int_{M}e^{-Cs}e{\big
(} f(y,s){\big )}dv_{g}(y),
$$
which implies
$$
e{\big (}f(x,s+t){\big )}\leq
Bt^{-\frac{Q}{2}}e^{Ct}\int_{M}e{\big (} f(y,s){\big )}dv_{g}(y).
\tag{5.31}
$$
By choosing a fixed $t=\varepsilon _{0}\in (0,\frac{\alpha
^{2}}{2})$, we get the estimate
$$
e{\big (}f(x,s+\varepsilon _{0}){\big )}\leq C(\varepsilon
_{0})E{\big (}f(\cdot ,s){\big )}
$$
where $C(\varepsilon _{0})$ is a constant depending on
$\varepsilon _{0}$. \qed
\enddemo

In view of Lemmas 5.3 and 5.4, one needs to estimate $E(f)$ for a
solution $f$ of the subelliptic harmonic map heat flow in order to
obtain a long-time existence result. Note that Proposition 4.1
implies
$$
\frac{d}{dt}E_{H}{\big (}f(\cdot ,t){\big )}=-\int_{M}\mid \tau
_{H} {\big (}f(\cdot ,t){\big )}\mid ^{2}dv_{g}\leq 0. \tag{5.32}
$$
Consequently $E_{H}{\big (}f(\cdot ,t){\big )}\leq E_{H}(\varphi
)$, where $\varphi $ is the initial map of $f$. Therefore it is
enough to estimate $E_{V}{\big (}f(\cdot ,t){\big )}$ for the
long-time existence.

\proclaim{Theorem 5.5} Let $(M,H,g_{H};g)$ be a compact
sub-Riemannian manifold and let $(N,h)$ be a compact Riemannian
manifold with nonpositive sectional curvature. Then for any map
$\varphi :M\rightarrow N$, the subelliptic harmonic map heat flow
(4.9) admits a global smooth solution $f:M\times \lbrack 0,\infty
)\rightarrow N$.
\endproclaim
\demo{Proof} By Schwarz inequality and the curvature assumption on
$N$, we get immediately from (3.29) that
$$
(\bigtriangleup _{H}-\frac{\partial }{\partial t})e_{V}(f)\geq
-C_{1}e_{H}(f)-C_{2}e_{V}(f)  \tag{5.33}
$$
for some positive constants $C_{1}$ and $C_{2}$. Integrating
(5.33) gives
$$
\aligned
\frac{d}{dt}E_{V}(f)&\leq CE_{H}(f)+C_{2}E_{V}(f) \\
&\leq CE_{H}(\varphi )+C_{2}E_{V}(f) \\
&=C_{1}+C_{2}E_{V}(f)
\endaligned
$$
which implies
$$
\int_{0}^{t}\frac{dE_{V}(f)}{C_{1}+C_{2}E_{V}(f)}\leq t.
$$
It follows that
$$
\ln \left( C_{1}+C_{2}E_{V}(f)\right) -\ln \left(
C_{1}+C_{2}E_{V}(\varphi )\right) \leq C_{2}t
$$
that is,
$$
E_{V}(f)\leq \frac{1}{C_{2}}\{e^{C_{2}t}\left(
C_{1}+C_{2}E_{V}(\varphi )\right) -C_{1}\}
$$
Hence we find that the solution $f(\cdot ,t)$ does not blow up at
any finite time. \qed
\enddemo
\vskip 0.4 true cm

{\bf 5.3 Eells-Sampson type results} \vskip 0.4 true cm

We will establish Eells-Sampson type results in following two
cases: the source manifolds are either step-$2$ sub-Riemannian
manifolds or step-$r$ sub-Riemannian manifolds whose
sub-Riemannian structures come from some Riemannian foliations.
\vskip 0.4 true cm {\bf 5.3.1 Step-2 sub-Riemannian manifolds}
\vskip 0.3 true cm

Recall that $T(\cdot,\cdot)$ denotes the torsion of the Bott
connection $\nabla ^{\frak{B}}$ on $(M,H,g_{H},g)$. Let $\pi
:S(V)\rightarrow M$ be the unit sphere bundle of the vertical
bundle $V$, that is, $S(V)=\{v\in V:$ $\parallel v\parallel
_{g}=1\}$. For any $v\in S(V)$, the $v$-component of $T(\cdot
,\cdot )$ is given by $T^{v}(\cdot ,\cdot )=\langle T(\cdot ,\cdot
),v\rangle $. Then we have a smooth function $\eta
(v)=\frac{1}{2}\Vert T^{v}\Vert _{g}^{2}:S(V)\rightarrow R$. Using
Lemma 3.1 and an adapted frame field $ \{e_{A}\}_{A=1,...,m+d}$
for $(M^{m+d},H,g_{H};g)$, we obtain
$$
\aligned \eta (v)&=\sum_{1\leq i<j\leq m}\left( T_{ij}^{\alpha
}\right) ^{2}\langle e_{\alpha },v\rangle ^{2}\\
&=\sum_{1\leq i<j\leq m}\langle \lbrack e_{i},e_{j}],v\rangle
^{2}. \endaligned\tag{5.34}
$$

\proclaim{Lemma 5.6} $H$ is $2$-step bracket generating if and
only if $\eta (v)>0$ for each $v\in S(V)$.
\endproclaim
\demo{Proof} For any $v\in S(V)$ with $\pi (v)=x$, we let $X,Y$ be
any local sections of $H$ around $x$. Writing $X=X^{i}e_{i}$ and
$Y=Y^{j}e_{j}$, we get
$$
\aligned \lbrack X,Y]_{x}&\equiv
X^{i}(x)Y^{j}(x)[e_{i},e_{j}]_{x}\phantom{b}
\mod\{H_{x}\} \\
&\equiv X^{i}(x)Y^{j}(x)\langle \lbrack e_{i},e_{j}],e_{\alpha
}\rangle _{x}e_{\alpha }(x)\phantom{b} \mod\{H_{x}\}
\endaligned
$$
Hence $H$ is $2$-step bracket generating for $TM$ if and only if
$$
span_{1\leq i,j\leq m}\{[e_{i},e_{j}]_{x}\}\equiv V_{x}\phantom{b
}\mod \{H_{x}\}
$$
at each point $x\in M$. By (5.34), this is equivalent to $\eta
(v)>0$. \qed
\enddemo

\proclaim{Lemma 5.7} Let $(M,H,g_{H};g)$ be a compact step-$2$
sub-Riemannian manifold and set $\eta _{\min }=\min_{v\in
S(V)}\eta (v)$. Let $N$ be a compact Riemannian manifold with
non-positive sectional curvature. Suppose $f:M\times \lbrack
0,\delta)\rightarrow N$ is a solution of the subelliptic harmonic
map heat flow. Then, for any given $t_{0}\in (0,\delta)$, we have
$$
E_{V}(f(\cdot ,t))\leq E_{V}(f(\cdot ,t_{0}))+\frac{4}{\eta _{\min
}}{\Big (}\int_{M}|\tau _{H}(f(\cdot ,t_{0}))|^{2}+CE_{H}{\big
(}f(\cdot ,t_{0}) {\big )}{\Big )}.
$$
for any $t\in (t_{0},\delta)$.
\endproclaim
\demo{Proof} The compactness of $M$ implies that $S(V)$ is
compact, so there exists a point $v\in S(V)$ such that $\eta
_{\min }=\eta (v)$. Since $H$ is $2$-step bracket generating, we
know from Lemma 5.6 that $\eta _{\min }>0$. Let $\varepsilon $ be
a fixed positive number with $\varepsilon \leq \frac{\eta _{\min
}}{4}$. From (3.31), (3.17) and (5.34), one has
$$
\aligned (\bigtriangleup _{H}-\frac{\partial }{\partial
t})e(f)&\geq -C_{\varepsilon }e_{H}(f)-\varepsilon
e_{V}(f)+(f_{ik}^{I})^{2}+\frac{1}{2}(f_{\alpha
k}^{I})^{2} \\
&\geq -C_{\varepsilon }e_{H}(f)-\varepsilon
e_{V}(f)+\sum_{I}\sum_{i<j}{\big (}
(f_{ij}^{I})^{2}+(f_{ji}^{I})^{2}{\big )}\\
&=-C_{\varepsilon }e_{H}(f)-\varepsilon
e_{V}(f)+\frac{1}{2}\sum_{I}\sum_{i<j}
{\big (} (f_{ij}^{I}+f_{ji}^{I})^{2}+(f_{ij}^{I}-f_{ji}^{I})^{2}{\big )} \\
&\geq -C_{\varepsilon }e_{H}(f)-\varepsilon e_{V}(f)+\frac{1}{2}
\sum_{I}\sum_{\alpha }\sum_{i<j}(f_{\alpha
}^{I})^{2}(T_{ij}^{\alpha })^{2}
\\
&=-C_{\varepsilon }e_{H}(f)-\varepsilon e_{V}(f)+\frac{1}{2}
\sum_{I}\sum_{\alpha }(f_{\alpha }^{I})^{2}\eta (e_{\alpha }) \\
&\geq -C_{\varepsilon }e_{H}(f)-\varepsilon
e_{V}(f)+\frac{1}{2}\eta _{\min }e_{V}(f).
\endaligned\tag{5.35}
$$
Integrating (5.35) over $M$ yields
$$
\aligned \frac{d}{dt}E(f)&\leq C_{\varepsilon
}E_{H}(f)+\varepsilon E_{V}(f)-\frac{
\eta _{\min }}{2}E_{V}(f) \\
&\leq C_{\varepsilon }E_{H}(f)-\frac{\eta _{\min }}{4}E_{V}(f).
\endaligned
$$
Consequently
$$
\frac{d}{dt}E_{H}(f)+\frac{d}{dt}E_{V}(f)+\frac{\zeta _{\min }}{4}
E_{V}(f)\leq C_{\varepsilon }E_{H}(f(\cdot ,t_{0})). \tag{5.36}
$$
By Corollary 4.6, we have
$$
\frac{d^{2}}{dt^{2}}E_{H}(f)\geq 0,
$$
which implies that
$$
\frac{d}{dt}E_{H}(f(\cdot ,t))\geq \frac{d}{dt}E_{H}(f(\cdot
,t_{0}))=-\int_{M}|\tau _{H}(f(\cdot ,t_{0}))|^{2}.  \tag{5.37}
$$
Set $A=\int_{M}|\tau _{H}(f(\cdot ,t_{0}))|^{2}+C_{\varepsilon
}E_{H}(f(\cdot ,t_{0}))$. From (5.36) and (5.37), it follows that
$$
\frac{d}{dt}E_{V}{\big (}f(\cdot ,t){\big )}+\frac{\eta _{\min
}}{4}E_{V}{\big (}f(\cdot ,t){\big )}\leq A
$$
that is,
$$
\frac{d}{dt}\left( e^{\frac{\eta _{\min }}{4}t}E_{V}{\big
(}f(\cdot ,t){\big )}\right) \leq Ae^{\frac{\eta _{\min }}{4}t}.
\tag{5.38}
$$
By integrating (5.38) over $[t_{0},t]$, we find
$$
e^{\frac{\eta _{\min }}{4}t}E_{V}{\big (}f(\cdot ,t){\big
)}-e^{\frac{\eta _{\min }}{4} t_{0}}E_{V}{\big (}f(\cdot
,t_{0}){\big )}\leq \frac{4A}{\eta _{\min }}{\big (} e^{\frac{\eta
_{\min }}{4}t}-e^{\frac{\eta _{\min }}{4}t_{0}}{\big )}.
$$
Hence
$$
\aligned E_{V}{\big (}f(\cdot ,t){\big )}&\leq e^{\frac{\eta
_{\min }}{4}(t_{0}-t)}E_{V}{\big (}f(\cdot
,t_{0}){\big )}+\frac{4A}{\eta _{\min }}(1-e^{\frac{\eta _{\min }}{4}(t_{0}-t)}) \\
&\leq E_{V}{\big (}f(\cdot ,t_{0}){\big )}+\frac{4A}{\eta _{\min
}}.
\endaligned
$$\qed
\enddemo

\proclaim{Theorem 5.8}Let $(M,H,g_{H};g)$ be a compact step-$2$
sub-Riemannian manifold and let $N$ be a compact Riemannian
manifold with non-positive sectional curvature. Then, for any
smooth map $\varphi :M\rightarrow N$, there exists a $C^{\infty }$
solution $f(x,t)$ of the subelliptic harmonic map heat flow (4.9)
on $M\times \lbrack 0,\infty )$. Moreover, there exists a sequence
$t_{i}\rightarrow \infty $ such that $f(x,t_{i})\rightarrow
f_{\infty }(x)$ uniformly, as $t_{i}\rightarrow \infty $, to a
$C^{\infty }$ subelliptic harmonic map $f_{\infty }:M\rightarrow
N$.
\endproclaim
\demo{Proof}Let $\frak{I}:N\hookrightarrow R^{K}$ be an isometric
embedding. Theorem 4.4 tells us that solving (4.9) is equivalent
to solving (4.19). In view of Theorem 5.1 and Lemmas 5.3, 5.4,
5.7, we conclude that (4.19) admits a global $C^{\infty \text{
}}$solution $u= \frak{I}\circ f:M\times \lbrack 0,\infty
)\rightarrow N\subset R^{K}$ with $f$ solving (4.9).

Now we investigate the convergence of $u$ as $t\rightarrow \infty
$. First, one observes that the compactness of $N$ and the uniform
boundedness of $e(u_{t})$ implies that the $1$-parameter family of
maps $u(\cdot ,t)$ form a uniformly bounded and equicontinuous
family of maps. Therefore, by Arzela-Ascoli Theorem, there exists
a sequence $t_{i}\rightarrow \infty $ such that
$$
u(\cdot ,t_{i})\rightarrow u_{\infty }(\cdot )  \tag{5.39}
$$
to a Lipschitz map $u_{\infty }:M\rightarrow N\subset R^{K}$.

Let us now deduce the equation which $\mid f_{t}\mid ^{2}=\mid
df(\frac{\partial }{\partial t})\mid ^{2}$ satisfies. By a direct
computation, using the commutation formulas (3.38) and (3.41), we
have
$$
\aligned \left( \bigtriangleup _{H}-\frac{\partial }{\partial
t}\right) (f_{t}^{I})^{2}&=2\left( f_{tk}^{I}\right)
^{2}+2f_{t}^{I}f_{tkk}^{I}-2\zeta
^{k}f_{t}^{I}f_{tk}^{I}-2f_{t}^{I}f_{tt}^{I} \\
&=2(f_{tk}^{I})^{2}+2f_{t}^{I}f_{ktk}^{I}-2\zeta
^{k}f_{t}^{I}f_{kt}^{I}-2f_{t}^{I}f_{tt}^{I} \\
&=2\left( f_{tk}^{I}\right) ^{2}+2f_{t}^{I}f_{kkt}^{I}-2f_{t}^{I}f_{k}^{K}%
\widetilde{R}_{KJL}^{I}f_{t}^{J}f_{k}^{L}-2\zeta
^{k}f_{t}^{I}f_{kt}^{I}-2f_{t}^{I}f_{tt}^{I} \\
&=2\left( f_{tk}^{I}\right) ^{2}+2f_{t}^{I}\left( f_{kk}^{I}-\zeta
^{k}f_{k}^{I}\right) _{t}-2f_{t}^{I}f_{k}^{K}\widetilde{R}
_{KJL}^{I}f_{t}^{J}f_{k}^{L}-2f_{t}^{I}f_{tt}^{I} \\
&=2\left( f_{tk}^{I}\right) ^{2}-2f_{t}^{I}f_{k}^{K}\widetilde{R}
_{KJL}^{I}f_{t}^{J}f_{k}^{L}.
\endaligned\tag{5.40}
$$
In terms of the curvature condition of $N$, (5.40) yields
$$
\left( \bigtriangleup _{H}-\frac{\partial }{\partial t}\right)
\mid f_{t}\mid ^{2}\geq 0.  \tag{5.41}
$$
By integrating (5.32) on any $[0,\delta]$, we get
$$
\int_{0}^{\delta}\int_{M}\mid f_{s}\mid ^{2}dv_{g}ds=E_H(\varphi
)-E_H(\delta)
$$
which implies that
$$
\int_{0}^{\infty }\int_{M}\mid f_{s}\mid ^{2}dv_{g}ds<\infty .
$$
Therefore there exists a sequence $s_{n}\rightarrow \infty $ such
that $\int_{M}\mid f_{s_{n}}\mid ^{2}dv_{g}\rightarrow 0$. From
Corollary 4.6, we see that
$$
\frac{d^{2}E_{H}(f_{t})}{dt^{2}}=-\frac{d}{dt}\{\int_{M}\mid
f_{t}\mid ^{2}dv_{g}\}\geq 0.
$$
Consequently $\int_{M}\mid f_{t}\mid ^{2}dv_{g}$ is decreasing in
$t$. Hence we find that
$$
\int_{M}\mid f_{t}\mid ^{2}dv_{g}\rightarrow 0  \tag{5.42}
$$
as $t\rightarrow \infty $. Clearly the function $\phi (x,t)=\mid
f_{s+t}\mid ^{2}$ also satisfies (5.41) for any given $s>0$.
Applying Lemma 2.4 to the function $\phi (x,t)$ for
$0<t<R_{0}^{2}$, we obtain
$$
\mid f_{s+t}\mid ^{2}\leq Bt^{-\frac{Q}{2}}\int_{M}\mid f_{s}\mid
^{2}dv_{g}. \tag{5.43}
$$
Then, for $t=\frac{R_{0}^{2}}{2}$, (5.43) gives that
$$
\mid f_{s+\frac{R_{0}^{2}}{2}}\mid ^{2}\leq
\frac{2^{\frac{Q}{2}}B}{R_{0}^{Q}}\int_{M}\mid f_{s}\mid
^{2}dv_{g}  \tag{5.44}
$$
for any $s>0$. From (5.42) and (5.44), it follows that
$$
\sup_{x\in M}\mid u_{t}\mid ^{2}(x,t)\rightarrow 0  \tag{5.45}
$$
as $t\rightarrow \infty $. Clearly (5.39) and (5.45) imply that
$u_{\infty }$ is a weak solution of (4.18). By Theorem 2.1, we can
now conclude that $u_{\infty }$ is smooth, that is, $f_{\infty }$
is a smooth subelliptic harmonic map from $M$ to $N$. \qed
\enddemo
\remark{Remark 5.1} It would be interesting to note that the
existence for Theorem 5.8 is independent of the choice of the
extension $g$ for $g_{H}$.
\endremark
\vskip 0.5 true cm {\bf 5.3.2 Riemannian foliations with basic mean
curvature vector} \vskip 0.3 true cm

Let $(M,H,g_{H};g)$ be a sub-Riemannian manifold corresponding to
a Riemannian foliation $\frak{F}$ on $(M,g)$ as in Example 1.4. A
foliation being Riemannian means that it is locally a Riemannian
submersion. In order to describe the local geometry of
$(M,g;\frak{F})$, we may assume temporarily that the foliation is
given by a Riemannian submersion $\pi :(M,g)\rightarrow
(Z,g_{Z})$. Then a vector field $X$ on $M$ is said to be
projectable if it is $\pi $-related to a vector field
$\widetilde{X}$ on $B$, that is, $ \widetilde{X}\circ \pi =\pi
_{\ast }(X)$.

\proclaim{Lemma 5.9} Let $(M,H,g_{H};g)$ be a sub-Riemannian
manifold corresponding to a Riemannian submersion $\pi
:(M,g)\rightarrow (Z,g_{Z})$. Let $X$ be a horizontal vector field
on $(M,H,g_{H};g)$. Then $X$ is projectable if and only if $\nabla
_{\xi }^{\frak{B}}X=0$ for any $\xi \in V$.
\endproclaim
\demo{Proof} Let $\Gamma(V)$ denote the space of vertical vector
fields. From [Mo], [GW], we know that a vector field $X$ on $M$ is
projectable if and only if $ [\xi ,X]\in \Gamma(V)$ for any $\xi
\in \Gamma(V)$, that is $\pi _{H}([\xi ,X])=0$. According to
(1.13), the lemma follows.\qed
\enddemo

In what follows, given a Riemannian submersion $\pi
:(M,g)\rightarrow (Z,g_{Z})$, a vector field $X$ on $M$ is said to
be basic if it is both horizontal and projectable.

\proclaim{Lemma 5.10} (cf. Lemma 1.4.1 in [GW]) Let
$(M,H,g_{H};g)$ be as in Lemma 5.9. If $X,Y\in \Gamma(M)$ are
basic, then so is $\nabla _{X}^{\frak{B}}Y$.
\endproclaim

Now we consider the general case that $(M,g;\frak{F})$ is a
Riemannian foliation. One says that $(M,g;\frak{F})$ is tense if
its mean curvature vector field $\zeta $ is parallel with respect
to $\nabla ^{\frak{B}}$ along the leaves, that is, $\nabla _{\xi
}^{\frak{B}}\zeta =0$ for any $ \xi \in V$. In view of Lemma 5.9,
we know that this condition means that $\zeta $ is (locally)
basic.

\proclaim{Lemma 5.11} Let $(M^{m+d},H,g_{H};g)$ be a compact
sub-Riemannian manifold corresponding a tense Riemannian foliation
$(M,g;\frak{F})$. Let $N$ be a compact Riemannian manifold with
non-positive sectional curvature. If $f:M\times \lbrack
0,\delta)\rightarrow N$ is a solution of the subelliptic harmonic
map heat flow, then $E_{V}(f_{t})$ is decreasing. In particular, $
E_{V}(f_{t})\leq E(\varphi )$.
\endproclaim
\demo{Proof} We first show that the curvature tensor of $\nabla
^{\frak{B}}$ satisfies
$$
\left. R_{j\alpha k}^{A}=0\right.  \tag{5.46}
$$
with respect to an adapted frame $\{e_{A}\}_{A=1}^{m+d}$. For any
point $ p\in M$, there exists a neighborhood $U$ of $p$ such that
the restriction of $\frak{F}$ to $U$ corresponds to a Riemannian
submersion $\pi :(U,g)\rightarrow (Z,g_{Z})$, since $\frak{F}$ is
Riemannian. Clearly we may choose an adapted frame field
$\{e_{A}\}_{A=1}^{m+d}$ such that $\{e_{1},...,e_{m}\}$ are basic
with respect to $\pi $, that is, $e_{j}\in \Gamma(U,H)$ and
$\nabla _{\xi }^{\frak{B}}e_{j}=0$ for any $\xi \in V$ ($1\leq
j\leq m$) due to Lemma 5.9. In view of Lemmas 5.9 and 5.10, we
also have $\nabla _{e_{\alpha }}^{\frak{B}}\nabla
_{e_{k}}^{\frak{B}}e_{j}=\nabla _{e_{\alpha }}^{\frak{B}}e_{j}=0$
and $\nabla _{\lbrack e_{\alpha },e_{k}]}^{\frak{B}}e_{k}=\nabla
_{\lbrack e_{\alpha },e_{k}]^{V}}^{\frak{B}}e_{k}=0$, where
$[e_{\alpha },e_{k}]^{V}$ denotes the vertical component of
$[e_{\alpha },e_{k}]$. Consequently
$$
\aligned
R_{j\alpha k}^{A}&=\langle R(e_{\alpha },e_{k})e_{j},e_{A}\rangle \\
&=\langle \nabla _{e_{\alpha }}^{\frak{B}}\nabla
_{e_{k}}^{\frak{B} }e_{j}-\nabla _{e_{k}}^{\frak{B}}\nabla
_{e_{\alpha }}^{\frak{B}}e_{j}-\nabla _{\lbrack e_{\alpha
},e_{k}]}^{\frak{B}}e_{k},e_{A}\rangle
\\
&=0. \endaligned \tag{5.47}
$$
In particular, one has $R_{k\alpha k}^{j}=0$. Using the
assumptions that $(M,g;\frak{F})$ is tense and $N$ has
non-positive curvature, we conclude from (3.29), (5.47) that
$$
\aligned (\bigtriangleup _{H}-\frac{\partial }{\partial
t})e_{V}(f_{t})&=(f_{\alpha k}^{I})^{2}+f_{\alpha }^{I}\zeta
_{,\alpha }^{k}f_{k}^{I}+f_{\alpha }^{I}f_{j}^{I}R_{k\alpha
k}^{j}-f_{\alpha }^{I}f_{k}^{K}\widehat{R}
_{KJL}^{I}f_{\alpha }^{J}f_{k}^{L} \\
&=(f_{\alpha k}^{I})^{2}-f_{\alpha }^{I}f_{k}^{K}\widehat{R}
_{KJL}^{I}f_{\alpha }^{J}f_{k}^{L} \\
&\geq 0. \endaligned \tag{5.48}
$$
Integrating (5.48) then gives this lemma. \qed
\enddemo
\remark{Remark 5.2} In [Dom], Dominguez showed that every
Riemannian foliation $\frak{F}$ on a compact manifold $M$ admits a
bundle-like metric $g$ for which the mean curvature vector field
$\zeta $ is basic. Hence tense Riemannian foliations exist in
abundance.
\endremark

Using Lemma 5.11 and a similar argument for Theorem 5.8, we obtain

\proclaim{Theorem 5.12} Let $(M,H,g_{H};g)$ be a compact
sub-Riemannian manifold corresponding to a tense Riemannian
foliation with the property that $H$ is bracket generating for
$TM$. Let $N$ be a compact Riemannian manifold with non-positive
sectional curvature. Then, for any smooth map $\varphi
:M\rightarrow N$, there exists a $C^{\infty }$ solution $f(x,t)$
of the subelliptic harmonic map heat flow (4.9) on $M\times
\lbrack 0,\infty )$. Moreover, there exists a sequence
$t_{i}\rightarrow \infty $ such that $f(x,t_{i})\rightarrow
f_{\infty }(x)$ uniformly, as $t_{i}\rightarrow \infty $, to a
$C^{\infty }$ subelliptic harmonic map $f_{\infty }:M\rightarrow
N$.
\endproclaim

Before ending this section, we would like to mention that Z.R.
Zhou [Zh2] announced an Eells-Sampson type result for subelliptic
harmonic maps from a sub-Riemannian manifold with vanishing
$\Gamma $-tensor. Here the $\Gamma $-tensor was introduced by
Strichartz in [St]. However, $\Gamma \equiv 0$ if and only if the
horizontal distribution $H$ is integrable.

\vskip 0.4 true cm

\heading{\bf 6. Hartman type Results}
\endheading
\vskip 0.3 true cm

First, we show the smoothness of a family of solutions to the
subelliptic harmonic map heat flow with a family of smooth maps as
its initial value. Our proof is similar to that in [Ha] for the
harmonic map heat flow and that in [RY] for the pseudo-harmonic
map heat flow, but with suitable modifications.

\proclaim{Lemma 6.1} Let $\varphi (x,\lambda ):M\times \lbrack
0,a]\rightarrow N\subset R^{K}$ be a smooth map and, for each
fixed $\lambda \in \lbrack 0,a]$, let $ u(x,t,\lambda)$ be a
solution of the subelliptic harmonic map heat flow on $M\times
[0,\delta)$ such that $u(x,0,\lambda)=\varphi (x,\lambda )$. Then
$u:M\times (0,\delta)\times (0,a)\rightarrow N$ is smooth.
\endproclaim
\demo{Proof} Suppose $u(x,t,\lambda)$ satisfies
$$
\cases \bigtriangleup _{H}u-\frac{
\partial u}{\partial t}=F(x,t,\lambda)\\
u(x,0,\lambda)=\varphi (x,\lambda )
\endcases\tag{6.1}
$$
for $(x,t,\lambda)\in M\times (0,\delta)\times (0,a)$, where
$F(x,t,\lambda)=\left( \Pi _{bc}^{a}\langle \nabla
^{H}u^{b},\nabla ^{H}u^{c}\rangle \right) $. First, we assert that
for any integer $l\ge 1$, $u(x,t,\lambda)$, $\partial
^{j}u/\partial \lambda ^{j}$ and $\nabla^{H}_x\partial
^{j}u/\partial \lambda ^{j}$ ($j=1,2,\cdots, l$) are continuous on
$M\times \lbrack 0,\delta)\times \lbrack 0,a]$. This can be proved
by a re-examination (and differentiations with respect to $\lambda
$) of the successive approximations used in the proof of the short
time existence theorem (Theorem 5.1). In terms of Theorem 2.1, we
see that for any fixed $\lambda$, $u(\cdot,\cdot,\lambda)\in
C^{\infty }(M\times (0,\delta),N)$ and all partial derivatives of
$u$ with respect to $(x,t)$ are bounded on any compact subsets of
$M\times (0,\delta)\times (0,a)$. Besides, by an inductive
argument on $l$ and the uniqueness theorem for the subelliptic
harmonic heat flow (Theorem 5.2), we see that $u(x,t,\lambda)$ is
smooth in $\lambda$ for each $(x,t)\in M\times (0,\delta)$, and
all partial derivatives of $u$ with respect to $\lambda$ are
bounded on any compact subsets of $M\times (0,\delta)\times (0,a)$
too. Therefore we may use the `joint smoothness lemma' in [RS]
(Lemma 6.2 on page 266 in [RS]) to conclude that $u:M\times
(0,\delta) \times (0,a)\rightarrow N$ is smooth. \qed
\enddemo

Next, we have the following lemma.

\proclaim{Lemma 6.2} Let $(M,H,g_{H};g)$ be a compact
sub-Riemannian manifold and $N$ be a compact Riemannain manifold
with non-positive sectional curvature. Let $ \varphi (x,\lambda
):M\times \lbrack 0,a]\rightarrow N$ be a family of smooth maps
and for fixed $\lambda $, let $f(x,\lambda ,t)$ be the solution of
the subelliptic harmonic map heat on $[0,\delta)$ such that
$f(x,0,\lambda )=\varphi (x,\lambda )$. Then for each $\lambda \in
\lbrack 0,a]$,
$$
\sup_{M\times \{t\}\times \{\lambda \}}\mid df(\frac{\partial
}{\partial \lambda })\mid ^{2}
$$
is non-increasing in $t$.
\endproclaim
\demo{Proof} For the map $f:M\times \lbrack 0,\delta)\times
\lbrack 0,a]\rightarrow N$, we define the following function
$$
\aligned Q(x,t,\lambda)&=\langle df(\frac{\partial }{\partial
\lambda }),df(\frac{
\partial }{\partial \lambda })\rangle \\
&=f_{\lambda }^{I}f_{\lambda }^{I}
\endaligned \tag{6.2}
$$
where $df(\frac{\partial }{\partial \lambda })=f_{\lambda
}^{I}\widetilde{e}_{I}$. In terms of (3.23), (3.38) and (3.41), we
deduce from (6.2) that
$$
\aligned \left( \bigtriangleup _{H}-\frac{\partial }{\partial
t}\right) Q&=2f_{\lambda k}^{I}f_{\lambda k}^{I}+2f_{\lambda
}^{I}f_{\lambda kk}^{I}-2\zeta
^{k}f_{\lambda }^{I}f_{\lambda k}^{I}-2f_{\lambda }^{I}f_{\lambda t}^{I} \\
&=2\left( f_{\lambda k}^{I}\right) ^{2}+2f_{\lambda
}^{I}f_{k\lambda k}^{I}-2\zeta ^{k}f_{\lambda }^{I}f_{\lambda
k}^{I}-2f_{\lambda
}^{I}f_{t\lambda }^{I} \\
&=2\left( f_{\lambda k}^{I}\right) ^{2}+2f_{\lambda
}^{I}f_{kk\lambda }^{I}-2f_{\lambda
}^{I}f_{k}^{K}\widetilde{R}_{KJL}^{I}f_{\lambda
}^{J}f_{k}^{L}-2\zeta ^{k}f_{\lambda }^{I}f_{k\lambda
}^{I}-2f_{\lambda
}^{I}f_{t\lambda }^{I} \\
&=2\left( f_{\lambda k}^{I}\right) ^{2}+2f_{\lambda
}^{I}(f_{kk}^{I}-\zeta ^{k}f_{k}^{I})_{\lambda }-2f_{\lambda
}^{I}f_{k}^{K}\widetilde{R}
_{KJL}^{I}f_{\lambda }^{J}f_{k}^{L}-2f_{\lambda }^{I}f_{t\lambda }^{I} \\
&=2\left( f_{\lambda k}^{I}\right) ^{2}+2f_{\lambda }^{I}\left(
f_{kk}^{I}-\zeta ^{k}f_{k}^{I}-f_{t}^{I}\right) _{\lambda
}-2f_{\lambda
}^{I}f_{k}^{K}\widetilde{R}_{KJL}^{I}f_{\lambda }^{J}f_{k}^{L} \\
&=2\left( f_{\lambda k}^{I}\right) ^{2}-2f_{\lambda }^{I}f_{k}^{K}\widetilde{R%
}_{KJL}^{I}f_{\lambda }^{J}f_{k}^{L} \\
&\geq 0.\endaligned \tag{6.3}
$$
Hence the maximum principle (Lemma 2.4) implies that if $0\leq
\tau \leq t<\delta$ , then
$$
\sup_{x\in M}Q(x,t,\lambda)\leq \sup_{x\in M}Q(x,\tau ,\lambda)
$$
for every fixed $\lambda \in \lbrack 0,a]$. Hence the desired
quantity is non-increasing. \qed
\enddemo

Suppose $f_{0}$ and $f_{1}$ are any two maps from $M$ to $N$. In
terms of the Riemannian distance $d_{N}$ of $N$, we have the
following distance between these two maps
$$
d_{N}^{\infty }(f_{0},f_{1})=\sup_{x\in M}d_{N}{\big
(}f_{0}(x),f_{1}(x) {\big )}.  \tag{6.4}
$$
Next, when $f_{0}$ and $f_{1}$ are homotopic, we may introduce the
homotopy distance between them as follows: If $F:M\times \lbrack
0,1]\rightarrow N$ is a smooth homotopy from $f_{0}$ to $f_{1}$,
so that $F(x,0)=f_{0}(x)$ and $ F(x,1)=f_{1}(x)$, then the length
of $F$ is defined by
$$
L(F)=\sup_{x\in M}\int_{0}^{1}\mid dF(\frac{\partial }{\partial
\lambda } )\mid _{(x,\lambda )}d\lambda .  \tag{6.5}
$$
One defines the homotopy distance $\widetilde{d}(f_{0},f_{1})$ to
be the infimum of the lengths over all homotopies from $f_{0}$ and
$f_{1}$. When $N$ has non-positively sectional curvature, the
homotopy distance can be attained by a smooth homotopy $G$ between
$f_{0}$ and $f_{1}$ in which $\lambda \mapsto G(x,\lambda )$ is a
geodesic for each $x\in M$, and in this case $L(G)=\sup_{x\in
M}\mid dG(\frac{\partial }{\partial \lambda })\mid $ for each
$\lambda \in \lbrack 0,1]$ (cf. [Jo2], [SY]). It is easy to see
that
$$
d_{N}^{\infty }(f_{0},f_{1})\leq \widetilde{d}(f_{0},f_{1}),
\tag{6.6}
$$
and if $d_{N}^{\infty }(f_{0},f_{1})<inj(N)$ (the injective radius
of $N)$, then $d_{N}^{\infty
}(f_{0},f_{1})=\widetilde{d}(f_{0},f_{1})$. Note that in order to
define $d_{N}^{\infty }(f_{0},f_{1})$ or $\widetilde{d}
(f_{0},f_{1}) $, we only need a Riemannian metric on $N$, while
$M$ can be any compact smooth manifold without any metric.

\proclaim{Proposition 6.3} Let $(M,H,g_{H};g)$ be a compact
sub-Riemannian manifold and let $N$ be a Riemannian manifold with
non-positive sectional curvature. Suppose $f_{0}(x,t)$ and
$f_{1}(x,t)$ are solutions of the subelliptic harmonic map heat
flow on $[0,\delta)$ with homotopic initial data. Then $t\mapsto
\widetilde{d}(f_{0}(\cdot ,t),f_{1}(\cdot ,t))$ is non-increasing.
\endproclaim
\demo{Proof} For any fixed $t_{0}\in \lbrack 0,\delta)$, let $F$
be the minimizing homotopy from $f_{0}(\cdot ,t_{0})$ to
$f_{1}(\cdot ,t_{0})$, that is, $L(F)= \widetilde{d}{\big
(}f_{0}(\cdot ,t_{0}),f_{1}(\cdot ,t_{0}){\big )}$. By Theorem
5.1, we have a solution $f(x ,t,\lambda)$ of the subelliptic
harmonic map heat flow on $[t_0,t_0+\delta_0 )$ for some $\delta_0
>0$ such that $f(x ,t_{0},\lambda)=F(x,\lambda )$. For any $t\in
\lbrack t_{0},t_{0}+\delta_0 )$, it is clear that $f(x,t,\lambda
)$ is a homotopy between $f_{0}(x,t)$ and $f_{1}(x,t)$. For any
$t\in \lbrack t_{0},t_{0}+\delta_0 )$, using Lemma 6.2, we derive
that
$$
\aligned \widetilde{d}{\big (}f_{0}(\cdot ,t),f_{1}(\cdot ,t){\big
)}&\leq L(f(\cdot ,t ,\cdot))=\sup_{x\in M}\int_{0}^{1}\mid
df(\frac{\partial }{
\partial \lambda })\mid _{(x ,t,\lambda)}d\lambda \\
&\leq \int_{0}^{1}\sup_{x\in M}\mid df(\frac{\partial }{\partial
\lambda })\mid _{(x ,t,\lambda)}d\lambda \\
&\leq \int_{0}^{1}\sup_{x\in M}\mid df(\frac{\partial }{\partial
\lambda }
)\mid _{(x,t_{0},\lambda)}d\lambda \\
&=\widetilde{d}{\big (}f_{0}(\cdot ,t_{0}),f_{1}(\cdot
,t_{0}){\big ).}
\endaligned\tag{6.7}
$$
This completes the proof of Proposition 6.3. \qed
\enddemo

\proclaim{Theorem 6.4} Let $(M,H,g_{H};g)$ be either as in Theorem
5.8 or Theorem 5.12. Suppose $(N,h)$ is a compact Riemannian
manifold with non-positive sectional curvature. Then the
subelliptic harmonic map heat flow (4.9) exists for all $t\in
\lbrack 0,\infty )$ and converges uniformly to a subelliptic
harmonic map $f_{\infty }$ as $t\rightarrow \infty $. In
particular, any map $\varphi \in C^{\infty }(M,N)$ is homotopic to
a subelliptic harmonic map.
\endproclaim
\demo{Proof} According to either Theorems 5.8 or 5.12, we know
that the subelliptic harmonic map heat flow (4.9) admits a global
solution $f:M\times \lbrack 0,\infty )\rightarrow N$, and there
exists a sequence $\{t_{k}\}$ such that $f(x,t_{k})$ converges
uniformly to a subelliptic harmonic map $f_{\infty }(x) $ as
$t_{k}\rightarrow \infty $.

The uniform convergence implies that $d_{N}^{\infty }{\big
(}f(\cdot ,t_{k}),f_{\infty }(\cdot ){\big )}<inj(N)$ for
sufficiently large $k$, and thus there is a unique minimizing
geodesic from $f(x,t_{k})$ to $f_{\infty }(x)$, which depends
smoothly on $x$. These geodesics define a homotopy from $f(\cdot
,t_{k})$ to $f_{\infty }(\cdot )$. This means that the maps
$f(\cdot ,t_{k})$ with large $k$ (and hence all, since $f(\cdot
,t)$ is continuous in $t$) are homotopic to $f_{\infty }$. In view
of Proposition 6.3, we have
$$
\widetilde{d}{\big (}f(\cdot ,t_{k}+t),f_{\infty }(\cdot ){\big
)}\leq \widetilde{d}{\big (}f(\cdot ,t_{k}),f_{\infty }(\cdot
){\big )}= d_{N}^{\infty }{\big (}f(\cdot ,t_{k}),f_{\infty
}(\cdot ){\big )}
$$
for all $t\geq 0$. Hence we conclude that the selection of the
subsequence is not necessary and that $f(\cdot ,t)$ uniformly
convergence to $f_{\infty } $ as $t\rightarrow \infty $. \qed
\enddemo

In previous existence results, the initial map $\varphi
:M\rightarrow N$ is assumed to be smooth. Similar to the case of
the harmonic map heat flow, we may take a continuous map as the
initial value for the subelliptic harmonic map heat flow.

\proclaim{Corollary 6.5} Let $M$ and $N$ be as in Theorem 6.4.
Then any continuous map $\varphi :M\rightarrow N$ is homotopic to
a subelliptic harmonic map $f$.
\endproclaim
\demo{Proof} One just need to smooth out the map $\varphi $ to a
smooth map $\widetilde{\varphi }$ such that $\widetilde{\varphi }$
is homotopic to $\varphi $ (cf. [Jo1], page 103-104). By applying
Theorem 6.4 to $\widetilde{\varphi }$, we get this corollary
immediately.\qed
\enddemo
\remark{Remark 6.1} Alternatively, one may check the proof for
local existence (Theorem 5.1), since after any positive time $t$,
the approximate solutions become automatically smooth. The
remaining arguments are as in Theorems 5.1, 5.8 and 5.12.
\endremark

\proclaim{Corollary 6.6} Let $M$ and $N$ be as in Theorem 6.4. Let
$\varphi :M\rightarrow N$ be a continuous map. Then the space of
subelliptic harmonic maps homotopic to $\varphi $ is connected,
and subelliptic harmonic maps in $[\varphi ]$ are all minimizers
of $E_{H}(\cdot )$ having the same horizontal energy.
\endproclaim
\demo{Proof} First, let us choose a minimizing sequence $\varphi
_{k}$ ($k=1,2,...$) in $[\varphi ]$ for $E_{H}(\cdot )$. Then we
get subelliptic harmonic maps $f_{k} $ ($k=1,2,...$) by the
preceding corollary. It follows from Lemmas 5.4, 5.7, 5.11 and
(5.32) that $e(f_{k})$ ($k=1,2,...$) are uniformly bounded. Hence
there exists a sequence of $\{f_{k}\}$ converges uniformly to a
Lipschitz map $f_{\min }$. Clearly $\frak{I}\circ f$ is a weak
solution of (4.18) with
$$E_{H}(f_{\min })=\lim_{t\rightarrow \infty} E_{H}(f_{k})$$ and
thus $f_{\min }$ is subelliptic by Theorem 2.1.

Now let $f$ be any subelliptic harmonic map in $[\varphi ]$. Then
there is a homotopy $F:M\times \lbrack 0,1]\rightarrow N$ between
$f$ and $f_{\min }$. It is known that $F$ determines a smooth
geodesic homotopy $G:M\times \lbrack 0,1]\rightarrow N$ between
these two maps. In [Zhou2], Zhou used the second variation formula
to show that each map in a geodesic homotopy between two
subelliptic harmonic maps has the same horizontal energy.
Consequently $E_{H}(G(\cdot ,t))=E_{H}(f_{\min })$. Therefore we
may conclude that each map $G(\cdot ,t)$ is a minimizing
subelliptic harmonic map for $E_{H}(\cdot )$, and the space of
subelliptic harmonic maps in $[\varphi ]$ is connected. \qed
\enddemo

From Examples 1.2, 1.3 and Theorems 5.8, 6.4, we immediately get

\proclaim{Corollary 6.7} Let $(M,H,g_{H})$ be either a compact
contact manifold or a compact quaternionic contact manifold with a
compatible metric $g$ and let $N$ be a compact Riemannian manifold
with non-positive sectional curvature. Then, for any continuous
$\varphi :M\rightarrow N$, there exists a $C^{\infty }$
subelliptic harmonic map $f_{\infty }:M\rightarrow N$ homotopic to
$\varphi $, which is a minimizer of $E_{H}$ in $[\varphi ]$.
\endproclaim

\remark{Remark 6.2} If $M$ is in particular a strictly
pseudoconvex CR manifold, the pseudoharmonic maps considered in
[ChC] and [RY], are subelliptic harmonic maps defined with respect
to the Webster metrics, while these metrics are only special
Riemannian extensions of the sub-Riemannian metrics determined by
the Levi forms. Hence, even in the CR case, the above Corollary
6.7 generalizes their results to the case that $g$ may be
arbitrary Riemannian extensions of the sub-Riemannian metrics (see
also Remark 5.1). This may provide some convenience for
considering further geometric analysis problems for subelliptic
harmonic maps on these manifolds.
\endremark

\vskip 0.5 true cm {\bf Acknowledgments}: The author would like to
thank Professor P. Cheng for helpful discussions.

\vskip 0.5 true cm 

\vskip 1 true cm

School of Mathematical Science

and

Laboratory of Mathematics for Nonlinear Science

Fudan University,

Shanghai 200433, P.R. China

\vskip 0.2 true cm yxdong\@fudan.edu.cn

\enddocument